\documentclass[12pt,oneside]{amsart}
\usepackage{amsmath,amsthm,amssymb}
\usepackage{color}
\usepackage{hyperref}
\usepackage{xcolor}
\usepackage{tikz}
\usepackage{comment}
\usepackage{indentfirst}
\usepackage{geometry}
\usepackage{mathrsfs}
\usepackage{tabularx}
\usepackage{multirow}
\usepackage{graphicx,color}
\usepackage{enumerate}
 \usepackage[numbers,sort&compress]{natbib}
\allowdisplaybreaks
\everymath{\displaystyle}

\theoremstyle{definition}
\newtheorem{lemma}{Lemma}[section]
\newtheorem{definition}[lemma]{Definition}
\newtheorem{proposition}[lemma]{Proposition}
\newtheorem{theorem}[lemma]{Theorem}
\newtheorem{corollary}[lemma]{Corollary}

\newtheorem{remark}[lemma]{Remark}

\newtheorem{conjecture}{Conjecture}
\numberwithin{equation}{section}
\allowdisplaybreaks

\makeatletter
\@namedef{subjclassname@2020}{%
  \textup{2020} Mathematics Subject Classification}
\makeatother

\title{On some new Ricci flow invariant curvature conditions}
\author{Yiyan Xu}
\address{School of Mathematics, Nanjing University, Nanjing, China, 210093}
\email{xuyiyan@nju.edu.cn}
\subjclass[2020]{53C20; 53E20}
\keywords{Ricci flow; Curvature operator; Sphere theorem; K\"ahler manifold}
\date{\today}

\begin{document}
\maketitle
\begin{abstract}
We would like to study new Ricci flow invariant curvature conditions. Specifically, we provide quantitative evidence for an unpublished conjecture  of  B\"ohm and  Wilking (see Conjecture \ref{BW-Conjecture}).    As an application, we study the topology of manifolds with pinched curvature.   
\end{abstract}

\section{Introduction and main results}
One of the fundamental  problems in Riemannian geometry is studying the interaction between   curvature and   topology of manifolds.  The famous sphere theorem states that a complete, simply connected, $1/4-$pinched manifold is homomorphic to a sphere. This result was proven by  Berger \cite{MR140054} and  Klingenberg \cite{MR139120} around 1960s using comparison techniques.  Moe recently,   Brendle and Schoen  \cite{MR2449060}  showed that the  homeomorphism can be sharped to diffeomorphism by using Ricci flow method,  and the theorem is referred  as the  Differentiable Sphere Theorem.

The Ricci flow was introduced by Hamilton \cite{MR0664497} in  1982.   In this fundamental paper  \cite{MR0664497}, Hamilton showed that on  three-dimensional closed manifolds,  the normalized Ricci flow  evolves the metrics with  strictly positive Ricci curvature
to  constant curvature limit metrics.   Consequently,  the manifolds are  diffeomorphic to  spherical space forms.   In the subsequent works of many authors, it turned out that  the Ricci flow   works for any dimension, provided that  the  initial curvature of the manifolds are  pinched in certain sense. Now the Ricci flow has been found to be a powerful tool to study this kind of the Differentiable Sphere Theorem.  

In \cite{MR0806701}, Huisken proved that $n$-dimensional closed manifolds are  diffeomorphic to  spherical space forms  if their curvature  do 
not deviate much from that of the spheres with constant sectional curvature, i.e., the norm of  Weyl curvature and the traceless Ricci tensor are small compared to the scalar curvature at each point.  In \cite{MR0862046} Hamilton showed that $4$-dimensional closed manifolds are  diffeomorphic to  spherical space forms  provided  that the   manifolds admit metrics with positive  curvature operators,   which was later extended by Chen \cite{MR1136125} to closed four-manifolds with two positive curvature operators. In a celebrated paper,  B\"ohm and  Wilking  \cite{MR2415394} further generalized Hamilton \cite{MR0862046} and Chen's  \cite{MR1136125} results to  manifolds with 2-positive
curvature operators in all dimensions.  Most importantly,  Hamilton \cite{MR0862046} introduced  the maximum principle for tensors and pinching sets for Ricci flow, which serve as a mandatory convergence criterion for Ricci flow today. 
While B\"ohm and  Wilking \cite{MR2415394} provided a generally effective method for  constructing   curvature pinching sets  by deforming curvature pinching families from known Ricci flow invariant curvature sets.   Soon thereafter,  many  other Ricci flow invariant curvature conditions  were studied,  such as
PIC1 and PIC2 \cite{MR2449060},  PIC \cite{MR3997128}, etc.  There is also a unified construction due to Wilking \cite{MR3065160} which recovers most of the Ricci flow invariant curvature conditions.

According to the  irreducible  decomposition theory \eqref{Irreducible-Curvature-Decomposition-1} for the space of Riemannian curvature operators, the  curvature operators  can be decomposed into three parts, 
\[R=  R_I+R_{{\rm Ric}_0}+W, \]
where $R_I$, $R_{{\rm Ric}_0}$, and $W$  
denote the scalar curvature part,   the traceless Ricci part, and the Weyl curvature part respectively, which will be defined precisely in 
\eqref{Curvature-Decomposition-1}.
In particular, $R=R_I$ corresponds to the curvature operator of the sphere $S^n$ with constant sectional curvature. 
 More precisely,    Huisken \cite{MR0806701} considered the following curvature cones:
\begin{equation*}
\Theta(\delta)=\Big\{R\in S_B^2(\mathfrak{so}(n))\,\Big| \|W\|^2+ \|R_{{\rm Ric}_0}\|^2 \leq  \delta \|R_I\|^2, {\rm scal}(R)>0\Big\},
\end{equation*} 
here $\|\cdot\|$ denotes the norm of curvature operators.  Note that  $\Theta(0)=\mathbb{R}^{+} I$, the cone of  curvature operators with constant positive sectional  curvature, and the cones $\Theta(\delta)$ are strictly convex.   Huisken   \cite{MR0806701} showed that  the curvature cones $\Theta(\delta)$ are preserved by Ricci  flow  for $\delta \leq \delta_n:=2/((n-2)(n+1))$ (in dimension 5, one need that $\delta_5= 1/10$). 
Furthermore, on  a compact Riemannian manifold $(M,g)$ with
the curvature operator at each point is  contained in the interior  of $\Theta(\delta)$  for some $\delta\in(0,\delta_n]$, Huisken   \cite{MR0806701} showed that  the normalized Ricci flow  evolves $g$ to a constant curvature limit metric (see also Margerin \cite{MR0840282}). Therefore,  $M$ is  diffeomorphic to a spherical space form.  
As pointed out  by Huisken   \cite{MR0806701}, the constant $\delta_n$, $n\geq 6$,  is the largest value of $\delta$ such that the interior of  
$\Theta(\delta)$ is contained in $\{R\in S_B^2(\mathfrak{so}(n))\, |\, R> 0\}$, the cone of positive curvature operators.  Therefore,  Huisken's   theorem in \cite{MR0806701}  was included in the latest  theorems in   \cite{MR0862046} and  \cite{MR2415394}.

An unpublished conjecture  of B\"ohm and Wilking says that:
\begin{conjecture}{\cite[Section 6.3.2]{Beitz-PhD-thesis}}\label{BW-Conjecture} Let $n\geq 12$, for $a\in[0,n/4]$,
the curvature cones 
\begin{equation}\label{B-W-Conjecture-Def-1}
\Omega(a)=\Big\{R\in S_B^2(\mathfrak{so}(n))\,\Big|\Big(\frac{n-2}{4}+a\Big)\|R\|^2\leq \frac{1}{4}|{\rm Ric}|^2,  {\rm scal} >0\Big\}.
\end{equation}
are invariant under the Hamilton ODE \eqref{Hamilton-ODE-Equ}  and  preserved by Ricci flow \eqref{Hamilon-Curvature-Evolution-1}.
\end{conjecture}
 With the    decomposition of Riemannian curvature operators \eqref{Irreducible-Curvature-Decomposition-1} again, the curvature cones $\Omega(a)$ in \eqref{B-W-Conjecture-Def-1} are also characterized by 
\begin{equation}\label{B-W-Cone-Def-1}
\begin{split}
\Omega(a) 
=\Big\{R\in S_B^2(\mathfrak{so}(n))\,\Big|a \|R_{{\rm Ric}_0}\|^2+\frac{n-2+4a}{4}\|W\|^2\leq \frac{n-4a}{4}\|R_I\|^2, {\rm scal} >0
\Big\}.
\end{split}
\end{equation}
Note that for $a > n/4$ 
the sets $\Omega(a)$  are empty. For $a =  n/4$, 
$\Omega(n/4)=\mathbb{R}^+ I$.
Moreover, for  $a\geq 0$  the cones 
$\Omega(a)$ are convex and in the case
that $a > 0$ even strictly convex. 

Beitz  considered the curvature cones $\Omega(a)$ in her  thesis \cite{Beitz-PhD-thesis} and derived some rigidity results for complete Einstein manifolds and complete shrinking gradient Ricci solitons based on the conjecture of B\"ohm and Wilking.  Some evidence for  the conjecture of B\"ohm and Wilking is the following \cite[Lemma 6.3.16]{Beitz-PhD-thesis}:   there exists some $\epsilon>0$,  but no estimate for the value of $\epsilon$, such that $\overline{\Omega(a)}=\Omega(a)\cup\{0\}$ respectively $\Omega(a)$ are invariant under the  Hamilton ODE  \eqref{Hamilton-ODE-Equ}  for $a\in[n/4-\epsilon,n/4]$.    

In this paper,  we provide  quantitative  evidence for  the conjecture of B\"ohm and Wilking by specifying an estimate for the value of $\epsilon$ which depends only on the dimension $n$ of the manifolds.

\begin{theorem}\label{Main-Theorem-1}There exists
\begin{equation}\label{BW-Conj-epsilon-n-Def}
\epsilon_n= \begin{cases}
     \frac{1}{n},&~\text{if}~n\geq 11, \\
    \frac{n^2(n-1)}{2(3n-2)(2n^2-4n+1)},&~\text{if}~4\leq n\leq 10,
\end{cases}
\end{equation}
such that the curvature cones $\overline{\Omega(a)}$ respectively $\Omega(a)$ are 
invariant under the Hamilton ODE  \eqref{Hamilton-ODE-Equ} and thus define Ricci flow invariant curvature conditions for $a\in[n/4-\epsilon_n,n/4]$. 
\end{theorem}

We would like to point out  that  Beitz's proof  relies on a second derivative test which    essentially   says that   the curvature operators with constant curvature  are     
  local attractors for   the Hamilton ODE \eqref{Hamilton-ODE-Equ}, while our proof relies on the  quantitative estimate  on algebraic curvature operators  which track back to Huisken \cite{MR0806701}.   
On the one hand, the curvature cones  $\Omega(a)$ for $a\in [n/4-\epsilon_n,n/4]$ in   Theorem \ref{Main-Theorem-1} are contained in the cone of curvature operators with positive Ricci curvature, see Corollary \ref{Ricci-Positive-Cond-1}. 
On the other hand,  it is easy to see that
\begin{equation}\label{Cur-Cone-Inc-Rel-1}
\Omega(n/4-\epsilon'_n)\subset \Theta(\delta_n)\subset \Omega(n/4-\epsilon_n),
\end{equation}
where $\epsilon'_n=1/(2(n-1))$, $\delta_n=2/((n-2)(n+1))$, and $\epsilon_n=1/n$.  In particular, when $n\geq 11$, the curvature cone  $\Omega(n/4-\epsilon_n)$  in our Theorem \ref{Main-Theorem-1} is slightly  larger than the largest curvature cone  $\Theta(\delta_n)$  in  Huisken's Theorem \cite{MR0806701}. However, the exact inclusion relationships between 
the curvature cone  $\Omega(n/4-\epsilon_n)$  in our Theorem \ref{Main-Theorem-1} and  the cone of positive curvature operators are not clear to us. 
%
%

Furthermore, with Hamilton \cite{MR0862046} and Bohm-Wilking's  \cite{MR2415394} terminology,  $\overline{\Omega(a)}_{a\in[n/4-\epsilon_n,n/4)}$ is a pinching family in the sense of Definition \ref{Pinching-family-Def}.   Therefore, with the  convergence theorem \ref{Hamilton-BW-CP-RF} for Ricci flow, 
we may conclude the following    the Differentiable Sphere Theorem.
\begin{theorem}\label{Main-Convergence-Thm} Let $n\geq 4$.
Suppose that $(M^n,g)$  is a compact Riemannian manifold such that
the curvature operator of $M$ at each point is   contained in   the interior of  $\Omega(n/4-\epsilon_n)$ (defined in \eqref{B-W-Cone-Def-1}), where $\epsilon_n$ is defined in \eqref{BW-Conj-epsilon-n-Def}.  Then the normalized Ricci flow  evolves $g$
to a constant curvature limit metric. Consequently,  $M$ is  diffeomorphic to a spherical space form. 
\end{theorem}

A K\"ahlerian analogue of the Sphere Theorem would be that if
$M$ is a compact K\"ahler manifold with positive bisectional curvature, then $M$ is biholomorphic to complex projective space $\mathbb{P}^m$; it is known as Frankel conjecture, proved by Mori \cite{MR0554387} and Siu-Yau \cite{MR0577360}. There has been much interest in obtaining a proof of
this using the K\"ahler Ricci flow.  On the one hand,   according to a result of Goldberg-Kobayashi \cite[Theorem 5]{MR0227901}, this amounts to show  the K\"ahler Ricci flow deforms the metric  with positive bisectional curvature to a K\"ahler Einstein metric. 
By assuming the existence of K\"ahler Einstein metric (the solution of the Frankel
conjecture),  Chen-Tian \cite{MR2219236}  showed that  the  normalized  K\"ahler Ricci flow evolves the metric with positive bisectional curvature to a metric with constant holomorphic sectional curvature (See also unpublished announcement of Perelman and Tian-Zhu \cite{MR2291916}).  On the other hand, although there is  a decomposition of a K\"ahler curvature operators  into  $U(m)$-irreducible subspaces, unlike the constant curvature operator of the round sphere, the K\"ahler curvature operator of the complex projective space with constant holomorphic sectional  curvature   is not a local attractor of the Hamilton ODE \eqref{Hamilton-ODE-Equ} restricted to the space of K\"ahler curvature operators (see \cite[Section 10]{MR1375255} or \cite[Section 6]{MR2415394}). Therefore,  the above result can not be proved   just carrying over the methods of real case to the K\"ahler case.

However,   we  can  show that  a K\"ahler manifold  $M$ of dimension $m$ has the same homotopy group with the complex projective space $\mathbb{P}^m$, provide its curvature  does not deviate much from that of $\mathbb{P}^m$ with constant holomorphic sectional curvature.  The   idea  can be  traced back  to   Kobayashi \cite{MR0154235}.  We know that a sphere $S^{2m+1}$ is a principal circle bundle over the complex projective space $\mathbb{P}^m$.  By generalizing
this situation, we can construct a principal circle bundle $P$ over $M$   such that  the Riemannian curvature operator  of   $P$  does not deviate much from that of  sphere $S^{2m+1}$. This allows us to run  Ricci flow  on $P$ and conclude that   $P$ is diffeomorphic with $S^{2m+1}$. The exact
homotopy sequences of the fiber bundle $\pi: P\rightarrow M$  will provide
the isomorphism  $\pi_i(M)=\pi_i(\mathbb{P}^m)$ for $i\geq 2$.

 In analogy to the Riemannian case,  the space of   K\"ahler curvature operators with  $U(m)$ action 
can also be  decomposed into three parts (see \eqref{Kah-Cur-irr-Dec-1}) and
\begin{equation}\label{Kahler-Curvature-irr-Dec-0}
K=K_E+K_{{\rm Ric}_0}+B,
\end{equation}
where $K_E, K_{{\rm Ric}_0}, B$ denote the scalar curvature part,   the traceless Ricci part and the Bochner-Weyl curvature part respectively,  which will be defined precisely in \eqref{Kah-Cur-Dec-For-1}.  In particular, $K=K_E$ corresponds to the curvature operator of the complex projective space $\mathbb{P}^m$ with constant holomorphic sectional curvature. 
 Additionally,  on   K\"ahler manifold $(M,g,J,\omega)$ of complex dimension $m$, we denote by $\bar{\lambda}={\rm scal}/(2m)$,   where ${\rm scal}$ denotes 
the Riemannian scalar curvature. The average of $\bar{\lambda}$ over $M$ is given by
\begin{equation}\label{Average-Scalar-Curvature-1}
\bar{\lambda}_0:=\frac{1}{{\rm Vol}(M)}\int_M\bar{\lambda}
=\frac{1}{m}\frac{ [c_1(M)]\wedge [\omega]^{n-1}}{[\omega]^n},
\end{equation}
 which is a constant depending on the K\"ahler class $[\omega]$ and the first Chern class $c_1(M)$.

\begin{theorem}\label{Kahler-Pinching-Bundle-3}  Let $m\geq 2$ and  $(M,g,J,\omega)$ be a compact  K\"ahler manifold of complex dimension $m$. Assume that the scalar curvature  is positive, i.e., $\bar{\lambda}>0$,
and  that the K\"ahler curvature $K$ at each point of $M$ satisfy the following  pinching condition,
\begin{equation}\label{Wilking-Kahler-Pinching-1}
 \begin{split}
 \Big(&\frac{3m-4}{8}+a\Big) \|K_{{\rm Ric}_0}\|^2+\Big(\frac{2m-1}{4}+a\Big)\|B\|^2\\
 &+\frac{m}{16(m+1)^2}\Big((12m^2-13)+4(6m+5)a +\frac{2 m+1 - 4 a }{2}\frac{1}{\epsilon}\Big) \Big(\bar{\lambda}-\bar{\lambda}_0\Big)^2\\
 &\leq \frac{2m+1-4a}{4} \frac{2m+1}{8(m+1)}\big(1-\frac{2}{2m+1}\epsilon\big) \|K_E\|^2,
  \end{split}
  \end{equation}
for some $0<\epsilon<(2m+1)/2$ and $(2m+1)/4-\epsilon_{2m+1}< a \leq (2m+1)/4$,   where $\epsilon_{2m+1}$ is defined in \eqref{BW-Conj-epsilon-n-Def}. Then   
  $\pi_i(M)\simeq \pi_i(\mathbb{P}^m)$ for  all $i$,  where  $\mathbb{P}^m$ denotes the complex projective space of complex dimension $m$.
\end{theorem}

\begin{remark}
Note that when $a$ is  close to $(2m+1)/4$ (e.g., $a>(2m+1)/4-\epsilon_{2m+1}'=(2m+1)/4-1/(4m)$, see \eqref{Cur-Cone-Inc-Rel-1}), the manifold $M$ has positive bisectional curvature, which implies that $M$ is biholomorphic to $\mathbb{P}^m$ by  Mori \cite{MR0554387} and Siu-Yau's \cite{MR0577360} solution of the Frankel conjecture. 

When the scalar curvature   is constant, i.e., $\bar{\lambda}\equiv\bar{\lambda}_0$,  the pinching condition  \eqref{Wilking-Kahler-Pinching-1} will be reduced to be a cleaner from (see the proof of Proposition \ref{PC-Kahler-Pinching-Hui-1}),
\begin{equation*}
\Big(\frac{3m-4}{8}+a\Big) \|K_{{\rm Ric}_0}\|^2+\Big(\frac{2m-1}{4}+a\Big)\|B\|^2 \leq \frac{2m+1-4a}{4} \frac{2m+1}{8(m+1)} \|K_E\|^2.
\end{equation*}
Furthermore, if the scalar curvature  of $M$ is constant, i.e., $\bar{\lambda}\equiv\bar{\lambda}_0$,  the Ricci 2-from is harmonic. If $M$ has positive  bisectional curvature,  then   $H^2(M;\mathbb{R})\simeq \mathbb{R}$ \cite{MR0227901}.  Thus the Ricci 2-from is proportional to the K\"ahler from of $M$, which means that $M$ is Einstein. Therefore,  $M$ is   isometric to $\mathbb{P}^m$ with the Fubini-Study metric by   a result of Goldberg-Kobayashi \cite[Theorem 5]{MR0227901}.  
\end{remark}

\noindent\textbf{Outline of the paper}:  We now give an outline of the paper.   In Section 2, we review    some relevant properties of algebraic curvature operators. In Section 3, we prove Theorems   \ref{Main-Theorem-1} and \ref{Main-Convergence-Thm} by  applying Hamilton's maximal principe for Ricci flow.  The key ingredient is  determining the  value of $a$ such that $Q(R)$ points to the   tangent cone of $\Omega(a)$ at $R\in \partial\Omega(a)$.    This requires  careful analysis and estimate of  the algebraic curvature operators. 
In Section 4, we  discuss the  principal circle bundles  over K\"ahler manifolds.
In Section 5, we study the topology of  K\"ahler manifolds with  pinched curvature operators by studying the principal circle bundles, 
and complete the proof  of Theorems  \ref{Kahler-Pinching-Bundle-3}. 

\section{Some properties for algebraic curvature operators}
\subsection{Algebraic curvature operators}\label{Al-Cur-Ope-1}
We will first introduce the space of algebraic curvature operators. Throughout the paper, we will adopt the notation from \cite{MR2415394} and \cite{MR3065160}.

 On a Riemannian manifold $(M^n,g)$,  the Riemannian curvature tensor is denoted by   \[R(x,y,z,w)=(R(x,y)z,w),~~x,y,z,w\in T_pM,\] here $R (x, y) =-\nabla_x\nabla_y+\nabla_y\nabla_x+\nabla_{[x,y]}.$ 

Then, for each $p\in M$,  one can choose an orthonormal basis $\{e_1,\cdots, e_n\}$ of $T_pM$ to build an isomorphism between $T_pM$
and $\mathbb{R}^n$, and also  isomorphism between $\Lambda^2T_pM$ and $\Lambda^2\mathbb{R}^n$.     We will identify $\Lambda^2\mathbb{R}^n$ with the Lie algebra $\mathfrak{so}(n)=\{X\in \mathbb{R}^{n\times n}: X+X^T=0\}$ by mapping the unit vector $e_i\wedge e_j$ onto the linear map $L(e_i\wedge e_j)$ of rank two which is a rotation with angle $\pi/2$ in the plane spanned by $e_i$ and $e_j$. We will identify $\alpha\wedge \beta$ with $L(\alpha\wedge \beta)$ in the following. Notice that under this identification, we have
\begin{equation}\label{Bi-vector-Inner-Product-1}
\begin{split}
\langle \alpha\wedge\beta,\gamma\wedge\delta\rangle_{\wedge^2V}&=\langle \alpha, \gamma\rangle\langle \beta, \delta\rangle-\langle \alpha, \delta\rangle\langle \beta, \gamma\rangle
=\langle \alpha\wedge\beta, \gamma\wedge\delta \rangle_{\mathfrak{so}(n)},
\end{split}
\end{equation}
thus $\Lambda^2V$ is isometric to $\mathfrak{so}(n)$ with the scalar product $\langle A, B\rangle =-1/2\,{\rm tr}(AB)$.  
 
Hereafter,  the Riemann curvature operator
 \begin{equation}\label{Cur-Ope-Rie-Ten-1}
{R}(e_i\wedge e_j)=\frac{1}{2}\sum_{k,l=1}^nR_{ijkl}e_k\wedge e_l,
\end{equation}
where $R_{ijkl}=R(e_i,e_j,e_k,e_l)$, 
is  viewed as a symmetric endomorphism of $\mathfrak{so}(n)$. 
We denote by $S_B^2(\mathfrak{so}(n))$  the vector space  of algebraic  curvature operators, that is the vector space of self-adjoint endomorphisms of $\mathfrak{so}(n)$ satisfying the Bianchi identity.   Moreover,  if ${R}\in S_B^2(\mathfrak{so}(n))$, 
we denote the curvature operator norm by  $\|R\|=\|R\|_{S_B^2(\mathfrak{so}(n))}=\langle {R}, {R}\rangle_{S_B^2(\mathfrak{so}(n))}^{1/2}$, where
\begin{equation}\label{Norm-Full-Curvature-Def-1}\begin{split}
\langle {R}, {R}\rangle_{S_B^2(\mathfrak{so}(n))}&:=\sum_{i<j}\langle {R}(e_i\wedge e_j), {R}(e_i\wedge e_j)\rangle=\frac{1}{4}\sum_{i,j,k,l=1}^n|R_{ijkl}|^2.
\end{split}
\end{equation}

The action of $O(n)$ on $\mathbb{R}^n$ induces the adjoint action of
$O(n)$ on   Lie algebra $\mathfrak{so}(n)$ and also the   action of $O(n)$ on $S_B^2(\mathfrak{so}(n))$, i.e., for $\sigma \in O(n)$, 
\begin{equation}\label{On-Action-Curvature-Def-1}
(\sigma R)(x\wedge y, z\wedge w)=R({\rm Ad}_{\sigma}(x\wedge y), {\rm Ad}_{\sigma}(z\wedge w))=R((\sigma x)\wedge (\sigma y), (\sigma z)\wedge (\sigma w)).
\end{equation}
For $n\geq 4$, the space of curvature operators  $S_B^2(\Lambda^2\mathbb{R}^n)$ has the following  orthogonal decomposition into irreducible  $O(n)$-invariant   subspaces (e.g., see \cite[Theorem 1.114]{MR2371700}):
\begin{equation}\label{Irreducible-Curvature-Decomposition-1}
S_B^2(\mathfrak{so}(n))=\langle I \rangle\oplus \langle {\rm Ric}_0 \rangle \oplus \langle W\rangle.
\end{equation}
Given  ${R}\in S^2_B(\mathfrak{so}(n))$, we let $R_{I}$,
$R_{{\rm Ric}_0}$ and $R_{W}$, denote the projections onto $\langle I \rangle$, $\langle {\rm Ric}_0 \rangle$ and $\langle W\rangle$, respectively.
Moreover,  let ${\rm Ric}$ denote the Ricci tensor of $R$,  ${\rm scal}={\rm tr}({\rm Ric})$ the scalar curvature
and  $W$ the Weyl tensor of $R$. 
Then
\begin{eqnarray}\label{Curvature-Decomposition-1}
    R_I=\frac{\bar{\lambda}}{n-1}{\rm id} \wedge {\rm id}, 
    \quad  
     R_{{\rm Ric}_0}=\frac{2}{n-2}{\rm Ric}_0 \wedge {\rm id}, \quad R_W=W, \label{RIRic}
\end{eqnarray}
where   $\bar{\lambda}={\rm scal}/n$,  ${\rm Ric}_0={\rm Ric}-\bar{\lambda}\,g$ is the traceless Ricci tensor,  and  the wedge product of two symmetric endomorphism $A, B: \mathbb{R}^n\rightarrow \mathbb{R}^n$ is defined by   (see \cite{MR2415394})
\begin{equation}\label{Wedge-Pro-Sym-1}
\begin{split}
A\wedge B:~&\Lambda^2\mathbb{R}^n\rightarrow \Lambda^2\mathbb{R}^n; \\
&x\wedge y\mapsto \frac{1}{2}\big(A(x)\wedge B(y)+B(x)\wedge A(y)\big).
\end{split}
\end{equation}
In particular,
\begin{equation}\label{Sym-End-Def-1}
\begin{split}
(A\wedge B)_{ijkl}=\frac{1}{2}(a_{ik}b_{jl}+b_{ik}a_{jl}-a_{il}b_{jk}-a_{jk}b_{il}).
\end{split}
\end{equation}
It is easy to see that  $A\wedge B=B\wedge A$ and $A\wedge B$ is a   curvature operator, while $I={\rm id} \wedge {\rm id}$ corresponds to the curvature operator of sphere $S^n$ scaled such that the sectional curvature is constant 1.

 Moreover, we have (e.g., compare with \cite[Lemma 2.3]{MR0806701})
\begin{equation}\label{Riemann-Ricci-Norm-1}
|{\rm Ric}|^2=|{\rm Ric}_0|^2+n\bar{\lambda}^2,
\end{equation}
\begin{equation}\label{curvature-norm-decomposition-3}
\|R_I\|^2=\|\frac{\bar{\lambda}}{n-1}{\rm id}\wedge {\rm id}\|^2=\frac{n}{2(n-1)}\bar{\lambda}^2,
\end{equation}
and 
\begin{equation}\label{curvature-norm-decomposition-5}
\|R_{{\rm Ric}_0}\|^2=\|\frac{2}{n-2}{\rm Ric}_0\wedge {\rm id}\|^2=\frac{1}{n-2}|{\rm Ric}_0|^2.
\end{equation}

Similarly,  let  $(M^m,g,J,\omega)$ be  a K\"ahler manifold of complex dimension $m$. For each $p\in M$, choose an orthonormal basis  $\{e_1,\cdots, e_{2m}\}$ of $T_pM$ of type $e_{m+i}=Je_i$, $i=1,\cdots, m$, set
\[E_i=\frac{1}{\sqrt{2}}(e_i-\sqrt{-1}Je_i), E_{\bar{i}}=\bar{E}_i=\frac{1}{\sqrt{2}}(e_i+\sqrt{-1}Je_i),~i=1,\cdots, m,\]
then $\{E_i\}_{i=1,\cdots,m}$ is an unitary basis of $T^{1,0}M$.  
We can identify $T_pM$
with  $\mathbb{R}^n$,  and $\mathbb{C}^m$,  $n=2m$.  The  action of complex structure 
 $J$   on $\Lambda^2\mathbb{R}^{n}$ gives the decomposition 
\[\Lambda^{2}\mathbb{R}^{n}=\Lambda^{1,1}\oplus \Lambda^{2,0}\oplus \Lambda^{0,2}.\]
If we identify $\Lambda^2(\mathbb{R}^{n})$ with $\mathfrak{so}(n)$,
then $\mathfrak{so}(n)\cap \Lambda^{1,1} \simeq \mathfrak{u}(m)$. 
 Note that the K\"ahler curvature operators $K$ satisfy $K(x,y)\circ J= J\circ K(x,y)$, then the  vector space   of all K\"ahler curvature operators  is then defined as the
space of elements of $S_B(\mathfrak{so}(n))$ which act trivially on the factor $\mathfrak{u}(m)^\perp$, can be identified with $S_B^2(\mathfrak{u}(m))$. For ${K}\in S_B^2(\mathfrak{u}(m))$,  we denote by ${\rm Ric}(K)$  the Ricci curvature of $K$, and $K_{i\bar{j}k\bar{l}}=K(E_i, E_{\bar{j}},E_k, E_{\bar{l}})$ and $K_{i\bar{j}}={\rm Ric}(K)(E_i, E_{\bar{j}})$. 

Given two hermitian endomorphisms $A=(a_{i\bar{j}}),B=(b_{i\bar{j}}): \mathbb{C}^m\rightarrow \mathbb{C}^m$,
we let $A\star B: \Lambda^{1,1}\rightarrow \Lambda^{1,1}$
denote the self-adjoint endomorphism of $\Lambda^{1,1}$ defined by (see \cite{MR3065160}\cite{MR322722})
\begin{equation}\label{Kahler-Wedge-Product-1}
A\star B_{i\bar{j}k\bar{l}}=-\frac{1}{2}(a_{i\bar{j}}b_{k\bar{l}}+a_{k\bar{l}}b_{i\bar{j}}+a_{i\bar{l}}b_{k\bar{j}}+a_{k\bar{j}}b_{i\bar{l}}), \quad 1\leq i,j,k,l\leq m,
\end{equation}
then $A\star B$ is a K\"ahler curvature operator. In particular,  $E={\rm id}\star {\rm id}$ with (e.g.,  see \cite{MR0227901})
\begin{equation}\label{Kahler-operator-CPn-1}
E_{ijkl}=\frac{1}{2}\big(\delta_{ik}\delta_{jl}-\delta_{il}\delta_{jk}+J_{ik}J_{jl}-J_{il}J_{jk}+2J_{ij}J_{kl}\big), \quad 1\leq i,j,k,l\leq 2m,
\end{equation}
corresponds to the curvature operator of  the complex projective pace $\mathbb{P}^m$ scaled such that the sectional curvature lies in the interval $[1/2,2]$. $E$ has the eigenvalue $m+1$ with multiplicity $1$,  the eigenvalue $1$ with multiplicity $m^2-1$ and the eigenvalue $0$ on $\mathfrak{u}(m)^\perp$ with multiplicity $m(m-1)$ (e.g., see section 5.2 in \cite{MR493867}).

As the analogue of the Riemannian case,  the space of K\"ahler curvature operators  $S_B^2(\mathfrak{u}(m))$ has the following  orthogonal decomposition into irreducible $U(m)$-invariant subspaces (\cite{MR0231313}\cite[Theorem 5.1]{MR322722}),
\begin{equation}\label{Kah-Cur-irr-Dec-1}
S_B^2(\mathfrak{u}(m))=\langle E \rangle\oplus \langle {\rm Ric}_0 \rangle \oplus \langle B\rangle.
\end{equation}
Given a K\"ahler curvature operator ${K}\in S^2_B(\mathfrak{u}(m))$, 
we let $K_{E}$,
$K_{{\rm Ric}_0}$ and $K_B$, denote the projections onto $\langle E \rangle$, $\langle {\rm Ric}_0 \rangle$ and $\langle B\rangle$, respectively. Note that  we denote  the projection of $K$ onto  $\langle {\rm Ric}_0 \rangle$ in   K\"ahler decomposition \eqref{Kah-Cur-irr-Dec-1} by $K_{{\rm Ric}_0}$ when it is unambiguous,  through it is not same as  the projection of $K$ onto  $\langle {\rm Ric}_0 \rangle$ in  Riemannian decomposition \eqref{Irreducible-Curvature-Decomposition-1}.   Then
\begin{equation}\label{Kah-Cur-Dec-For-1}
K_E=\frac{\bar{\lambda}}{m+1}\, {\rm id}\star {\rm id}, \quad K_{{\rm Ric}_0}=\frac{2}{m+2}\, {\rm Ric}(K)_0 \star {\rm id},\quad K_B=B,
\end{equation}
here  $\bar{\lambda}=\frac{1}{2m}{\rm scal}$,   ${\rm scal}$ denotes 
the Riemannian scalar curvature, ${\rm Ric}(K)_0={\rm Ric}(K)-\bar{\lambda} g$ denotes  the traceless Ricci curvature, $B$ denotes  the Bochner-Weyl curvature.

If ${K}\in S_B^2(\mathfrak{u}(m))$,  set
\[\begin{split}
\langle {K}, {K}\rangle_{S_B^2(\mathfrak{u}(m))}&:=\sum_{i,j,k,l=1}^mK_{i\bar{j}k\bar{l}}\overline{K_{i\bar{j}k\bar{l}}},
\end{split}\]
then a straightforward computation shows that (e.g., see \cite[Proposition  3.4]{MR3856847})
\begin{equation}\label{Kahler-Riemann-Norm-1}
\begin{split}
\sum_{i,j,k,l=1}^mK_{i\bar{j}k\bar{l}}\overline{K_{i\bar{j}k\bar{l}}}=\frac{1}{4}\sum_{i,j,k,l=1}^{2m} K_{ijkl}^2, 
\end{split}
\end{equation}
i.e, $\|K\|_{S_B^2(\mathfrak{u}(m))}=\|K\|_{S_B^2(\mathfrak{so}(2m))}$. Thus,  we will still  denote the (K\"ahler) curvature operator norm by  $\|K\|$  as in the real case.
Moreover, we have 
  \begin{equation}\label{Ricci-Norm-RK-Rel-1}
  \begin{split}\sum_{i,j=1}^mK_{i\bar{j}}\overline{K_{i\bar{j}}}
=\frac{1}{2}\sum_{i,j=1}^{2m} {\rm Ric}(K)_{ij}^2.
   \end{split}
   \end{equation}
We  denote by  $|{\rm Ric}(K)|_K:=\big(\sum_{i,j=1}^mK_{i\bar{j}}\overline{K_{i\bar{j}}}\big)^{\frac{1}{2}}$, then we have  $|{\rm Ric}(K)|_K=1/\sqrt{2}\,|{\rm Ric}(K)|$ and
 
\begin{equation}\label{Kaehler-Ric-Norm-1}
|{\rm Ric}|^2=|{\rm Ric}_0|^2+2m\bar{\lambda}^2,\quad |{\rm Ric}|_K^2=|{\rm Ric}_0|_K^2+m\bar{\lambda}^2.
\end{equation}
Equipped with the above notation,  we have  that 
 \begin{equation}\label{Kaehler-Curvature-Scalar-Norm-1} 
 \begin{split}
\|K_{E}\|^2&=\|\frac{\bar{\lambda}}{m+1}\, {\rm id}\star {\rm id}\|^2 = \frac{2m}{m+1}\bar{\lambda}^2,
\end{split}
\end{equation}
and
\begin{equation}\label{Kaehler-Curvature-Ric-Norm-1}
\begin{split}
\|K_{{\rm Ric}_0}\|^2&= \|\frac{2}{m+2}{\rm Ric}(K)_0 \star {\rm id}\|^2=\frac{4}{m+2} |{\rm Ric}(K)_0|_K^2.
\end{split}
\end{equation}

\subsection{Curvature pinching along Ricci flow}\label{RF-Inv-Cur-Cone-1}
On a  Riemannian manifold $(M,g_0)$, the Ricci flow is the geometric evolution equation
\begin{equation}\label{Ricci-Flow-Evo-Equ-1}
\frac{\partial g}{\partial t}=-2{\rm Ric}(g),
\end{equation}
starting at $g_0$. Using moving frames, 
the evolution equation of curvature operators $R=R_{g(t)}$ of $(M,g(t))$ along the Ricci flow \eqref{Ricci-Flow-Evo-Equ-1} is (\cite{MR0862046})

\begin{equation}\label{Hamilon-Curvature-Evolution-1}
 \frac{\partial}{\partial t} R=\Delta R +2Q(R),
\end{equation}
where  $Q$ is the quadratic vector field on $S_B^2(\mathfrak{so}(n))$  defined by
\[Q(R)=R^2+R^\sharp.\]
Here,  $R^2$ is  the square of $R\in S^2(\mathfrak{so}(n))$ seen as an endomorphism of $\mathfrak{so}(n)$, while
$R^\sharp=R\sharp R$ is defined in the following way,
\[\begin{split}
\sharp:\quad & S^2(\mathfrak{so}(n))\times S^2(\mathfrak{so}(n))\rightarrow S^2(\mathfrak{so}(n))\\
&({R},{S})\mapsto {R}\sharp {S}={\rm ad}\circ ({R}\wedge {S})\circ {\rm ad}^*,
\end{split} \]
where ${\rm ad}: \Lambda^2\mathfrak{so}(n)\rightarrow \mathfrak{so}(n), u\wedge v\mapsto [u,v]$ denotes the adjoint representation of $\mathfrak{so}(n)$ and ${\rm ad}^*$ is its dual.  We would like to mention that ${R}\sharp {S}={S}\sharp {R}$, and \begin{equation}\label{Ham-Sharp-Def-3}
\begin{split}
\langle {R}\sharp {S}(\varphi),\varphi\rangle
&=-\frac{1}{2}{\rm tr}({\rm ad}_\varphi \circ {R}\circ {\rm ad}_\varphi\circ {S}).
\end{split}
\end{equation}
In particular, with respect to normal frame,  we have 
\begin{align}
R^2_{ijkl}&=({R}\circ {R})_{ijkl}=\frac{1}{2}R_{ijpq}R_{klpq},  \label{Cur-Ope-Squ-Loc-For-1}\\
{R}^\sharp_{ijkl}&={R}\sharp {R}_{ijkl}=R_{ipkq}R_{jplq}-R_{iplq}R_{jpkq}.\label{Cur-Ope-Sharp-Loc-For-1}
\end{align}
Moreover, we have 
(see \cite[(6)]{MR2415394})
\begin{align}\label{Ricci-QR-For-1}
{\rm Ric}(Q(R))_{ij}
   &=\sum_{k=1}^n(R^2+R^\sharp)_{ikjk}
=\sum_{p,q=1}^nR_{ipjq}R_{pq}=\mathring{R}({\rm Ric})_{ij},
\end{align}
and
\begin{equation}\label{Cur-Ope-Tri-Loc-For-1}
\begin{split}
\langle Q(R), R\rangle&=\langle {R}^2+{R}^\sharp, {R}\rangle=\frac{1}{8}(R_{ijpq}R_{klpq}R_{ijkl}+4R_{ipkq}R_{jplq}R_{ijkl}).
\end{split}
\end{equation}

Now we turn to  Hamilton's maximum principle for Ricci flow. We consider the closed convex curvature set $C\subset S_B^2(\mathfrak{so}(n))$ which is invariant under the action $O(n)$ defined in \eqref{On-Action-Curvature-Def-1}.  Thanks to the $O(n)$-invariance of $C$, we can say the curvature operator of a manifold $(M,g)$ at each point $p$ is  contained in $C$ under the isomorphism between $T_pM$ and $\mathbb{R}^n$ that we build  at the beginning of Section \ref{Al-Cur-Ope-1}, which is independent of the basis of $T_pM$.

\begin{proposition}[\cite{MR0862046}]\label{Hamilton'-maximum-principle}
A closed convex $O(n)$-invariant subset $C\subset S_B^2(\mathfrak{so}(n))$ which is invariant under the ordinary differential equation 
\begin{equation}\label{Hamilton-ODE-Equ}
\frac{d}{dt}R=Q(R),
\end{equation}
that is, the  solutions of  ODE \eqref{Hamilton-ODE-Equ}  which start inside   $C$  stay in $C$,  
defines a Ricci flow invariant curvature condition, that is, the Ricci flow evolves metrics whose curvature operator at each point   is contained in $C$ into metrics with the same property. 
\end{proposition}
The following property is somewhat a criterion for a  Ricci flow invariant curvature set in terms of its tangent cone. 
\begin{proposition}\label{RF-Inv-Tangent-Cone-1}
A closed convex $O(n)$-invariant subset $C\subset S_B^2(\mathfrak{so}(n))$   is invariant under the Hamilton ODE  \eqref{Hamilton-ODE-Equ} if and only if for all $R\in \partial C$ we have $Q(R)\in T_RC$, the tangent cone of $C$ at $R$.
\end{proposition}

Hamilton \cite{MR0862046} established a general convergence criterion for Ricci flow, in which the  mandatory ingredients  are the construction of a so-called curvature pinching sets. B\"ohm and Wilking \cite{MR2415394} provided an effective method for  constructing   curvature pinching sets  from curvature pinching families. The idea can be tracked back to Hamilton  \cite{MR0664497}, although it was not really stressed. The following  definition of   pinching family comes from B\"ohm and Wilking \cite{MR2415394}.

\begin{definition}[Pinching family \cite{MR2415394}]\label{Pinching-family-Def}
We call a   family of curvature  cones  $C(s)_{s\in [0,1)}\subset
S_B^2(\mathfrak{so}(n))$ a pinching family, if
\begin{enumerate}[(1)]
\item  $C(s)$ is closed convex $O(n)$-invariant cones of full
dimension for all $s\in [0,1)$,
\item $s\mapsto C(s)$ is continuous;
\item each $R\in C(s)\setminus \{0\}$ has positive scalar curvature,
\item $C(s)$ converges in the pointed Hausdorff topology to the one-dimensional cone $\mathbb{R}^+ I$ as $s\rightarrow 1$.
\item $Q(R)=R^2+R^\sharp$ is contained in the interior of the tangent cone of $C(s)$ at $R$ for all $R \in C(s)\setminus \{0\}$ and all $s\in (0,1)$.
\end{enumerate}
\end{definition}

\begin{theorem}\cite[Theorem 5.1]{MR2415394}\label{Hamilton-BW-CP-RF}
Let $C(s)_{s\in[0,1)}\subset S_B^2(\mathfrak{so}(n))$  be a pinching family of closed convex cones,  $n\geq 3$. 
Suppose that $(M,g)$  is a compact Riemannian manifold such that
the curvature operator of $M$  at each point is  contained in  the interior of $C(0)$.  Then the normalized Ricci flow  evolves $g$
to a constant curvature limit metric.
\end{theorem}

\subsection{Some algebraic identities for curvature operators}
\begin{lemma}\cite[Lemma 2.1]{MR2415394}\label{BW-Ricci-Identity-1} Let $R \in S^2_B(\mathfrak{so}(n))$, then
\begin{equation}\label{BW-Identity-1}
R+R\sharp  I={\rm Ric}(R)\wedge {\rm id}.
\end{equation}
\end{lemma}

We say that a curvature operator $R$ is of Ricci type,
if $R=R_I+R_{{\rm Ric}_0}$, i.e. $R_W=0$.
\begin{lemma}\cite[Lemma 2.2]{MR2415394} \label{BW-Riccitype-Cur-1}
Let $R \in S^2_B(\mathfrak{so}(n))$ be a curvature operator of Ricci type,
then
\begin{equation}\label{BW-Q-Identity-2}
\begin{split}
  R^2+R^\sharp&=\frac{1}{n-2}{\rm Ric}_0\wedge {\rm Ric}_0
              +\frac{2\bar{\lambda}}{n-1}{\rm Ric}_0\wedge {\rm id}
         -\frac{2}{(n-2)^2}({\rm Ric}_0^2)_0\wedge {\rm id}\\
            &+\frac{\bar{\lambda}^2}{n-1}I
           +\frac{|{\rm Ric}_0|^2}{n(n-2)}\, I\,.
\end{split}
\end{equation}
Moreover
\begin{eqnarray}
  \bigl(R^2+R^\sharp\bigr)_{W}
   &=&
     \frac{1}{n-2}\bigl({\rm Ric}_0 \wedge {\rm Ric}_0\bigr)_{W}, \label{BW-Q-W-Identity-2}\\
  {\rm Ric}( R^2+R^\sharp)
   &=&
     -\frac{2}{n-2} ({\rm Ric}_0^2)_0+\frac{n-2}{n-1}\bar{\lambda} {\rm Ric}_0
     + \bar{\lambda}^2 {\rm id} +\frac{|{\rm Ric}_0|^2}{n}{\rm id}\,. \label{BW-Q-Ric-Identity-2}
 \end{eqnarray}
\end{lemma}

We will   denote the bilinear map by the same notation $Q$:  \[Q(R,S):=\frac{1}{2}\big(Q(R+S)-Q(R)-Q(S)\big)=\frac{1}{2}(RS+SR)+R\sharp S, ~{R},{S} \in S^2(\mathfrak{so}(n)).\] 
Recall that the trilinear map (see \cite[(5)]{MR2415394}),
\begin{equation}
{\rm tri}({R},{S},{T})={\rm tr}\big({R}{S}+{S}{R}+2{R}\sharp {S}){T}\big)=2\langle Q(R,S), T\rangle,~{R},{S},{T}\in S^2(\mathfrak{so}(n)),
\end{equation}
  is symmetric   in all three variables. 

The following lemma was essentially proved in \cite{MR2415394} (see also  \cite[Lemma 4.3]{MR2443993}).

\begin{proposition}\label{Q-R-Kahler-Decomposition-1} The evaluation of  bilinear map   $Q$    on different parts of the decomposition of algebraic curvature operators  have the following properties: 
\begin{enumerate}[(1)]
  \item If $R,S\in \langle I\rangle$, then $Q(R,S)\in  \langle I\rangle$;  
  \item If $R\in \langle I\rangle$, $S\in \langle {\rm Ric}_0\rangle$,  then $Q(R,S)\in \langle {\rm Ric}_0\rangle$;
   \item If $R\in \langle I\rangle$, $S\in \langle W\rangle$,  then $Q(R,S)=0$;
 \item If $R, S\in \langle W\rangle$, then $Q(R,S)\in \langle W\rangle$;
  \item If $R\in \langle {\rm Ric}_0\rangle, S\in \langle W\rangle$, then $Q(R,S)\in \langle {\rm Ric}_0\rangle$.
\end{enumerate}

\end{proposition}

\begin{proof} By \eqref{BW-Identity-1},  if $R,S\in \langle I\rangle$, then $Q(R,S)\in  \langle I\rangle$;  If $R\in \langle I\rangle$, $S\in \langle {\rm Ric}_0\rangle$,  then $Q(R,S)\in \langle {\rm Ric}_0\rangle$; If $R\in \langle I\rangle$, $S\in \langle W\rangle$,  then $Q(R,S)=0$. The conclusions  $(1), (2), (3)$ are checked.

 If $R\in \langle W\rangle$, by \eqref{Ricci-QR-For-1}, we have 
${\rm Ric}(Q(R))=\mathring{R}({\rm Ric})=0$, ${\rm scal}(Q(R))=|{\rm Ric}|^2=0$, thus $Q(R)=R^2+R^\sharp\in \langle W\rangle$; Consequently, If $R, S\in \langle W\rangle$, then $Q(R,S)=\frac{1}{2}(Q(R+S,R+S)-Q(R)-Q(S))\in \langle W\rangle$, which proves $(4)$.

If $R\in \langle {\rm Ric}_0\rangle, S\in \langle W\rangle$,  by applying then the symmetric property of the trilinear map ${\rm tri}$, we obtain from $(1), (2), (3), (4)$ that 
\[\langle Q(R,S), I\rangle=\frac{1}{2}{\rm tri}(R,S,I)=\langle Q(R,I), S\rangle=0,\]
and
\[\langle Q(R,S), T\rangle=\frac{1}{2}{\rm tri}(R,S,T)=\langle Q(S,T), R\rangle=0, \forall T\, \in \langle W\rangle,\]
therefore $Q(R,S)\in \langle {\rm Ric}_0\rangle$, and conclusion $(5)$ is verified.  
\end{proof}

In our paper, we will use the following identity, which was derived  by Huisken in 
his proof of the Theorem 3.3 in \cite{MR0806701}. Since this identity was not written down as an independent lemma in \cite{MR0806701},   for readers' convenience, we include  calculations here.
\begin{lemma}  For $R\in S^2(\mathfrak{so}(n))$, we have 
\begin{equation}\label{Tri-R-Full-Huisken-1}
\begin{split}
\langle Q(R), {R}\rangle&=\frac{n}{2(n-1)}\bar{\lambda}^3+\frac{3}{2(n-1)}\bar{\lambda}|{\rm Ric}_0|^2+\frac{3}{n-2}\langle{\rm Ric}_0\wedge {\rm Ric}_0,W\rangle\\
&\quad -\frac{2}{(n-2)^2}\langle  ({\rm Ric}_0^2)_0, {\rm Ric}_0\rangle+\langle Q(W), {W}\rangle.
\end{split}
\end{equation}
\end{lemma}
\begin{proof} First,  by applying   the symmetric property of the trilinear map ${\rm tri}$ and also  \eqref{BW-Ricci-Identity-1}, we have 
\begin{equation}\label{Tri-R-Full-Cal-Aux-1}
\langle Q(R),R_I\rangle  =\frac{1}{2}{\rm tri}(R,R_I,R) 
 =\frac{\bar{\lambda}}{n-1}\langle R+R\sharp I, R\rangle=\frac{\bar{\lambda}}{n-1}\langle {\rm Ric}(R)\wedge {\rm I}, R\rangle.
 \end{equation}
According to the curvature decomposition \eqref{Curvature-Decomposition-1},  we derive from \eqref{Tri-R-Full-Cal-Aux-1} to have 
\begin{equation}\label{Tri-R-Full-Cal-Aux-3}
\begin{split}
\langle Q(R),R_I\rangle 
&=\frac{\bar{\lambda}}{n-1}\langle  \bar{\lambda}\,{\rm id} \wedge {\rm id}+{\rm Ric}_0\wedge {\rm I},
\frac{\bar{\lambda}}{n-1}{\rm id} \wedge {\rm id}
+\frac{2}{n-2}{\rm Ric}_0 \wedge {\rm id}+W\rangle\\
&=\frac{\bar{\lambda}}{2(n-1)}\big(|{\rm Ric}_0|^2+n\bar{\lambda}^2\big).
\end{split}
\end{equation}
Next, with Lemma \ref{Q-R-Kahler-Decomposition-1}, we conclude that 
\begin{equation}\label{Tri-R-Full-Cal-Aux-5}
\begin{split}
\langle Q(R)&,R_{{\rm Ric}_0}\rangle=\frac{1}{2}{\rm tri}(R,R,R_{{\rm Ric}_0})\\
&=\frac{1}{2}{\rm tri}(R_I,R_{{\rm Ric}_0},R) +\frac{1}{2}{\rm tri}(R_{{\rm Ric}_0},R_{{\rm Ric}_0},R)+\frac{1}{2}{\rm tri}(W,R_{{\rm Ric}_0},R)\\
&={\rm tri}(R_I,R_{{\rm Ric}_0},R_{{\rm Ric}_0})+{\rm tri}(R_{{\rm Ric}_0},R_{{\rm Ric}_0},W)+\frac{1}{2}{\rm tri}(R_{{\rm Ric}_0},R_{{\rm Ric}_0},R_{{\rm Ric}_0})\\
&=2\langle Q(R_{{\rm Ric}_0}), R_I \rangle +2\langle Q(R_{{\rm Ric}_0}), W \rangle+\langle Q(R_{{\rm Ric}_0}), R_{{\rm Ric}_0} \rangle.
\end{split}
\end{equation}
By applying \eqref{BW-Q-Identity-2}  and \eqref{BW-Q-Ric-Identity-2} in Lemma \ref{BW-Riccitype-Cur-1} to $R_{{\rm Ric}_0}$ respectively,  we have 
\[Q(R_{{\rm Ric}_0})=\frac{1}{n-2}{\rm Ric}_0\wedge {\rm Ric}_0
         -\frac{2}{(n-2)^2}({\rm Ric}_0^2)_0\wedge {\rm id}           +\frac{|{\rm Ric}_0|^2}{n(n-2)}\, I,\]
and \[
{\rm Ric}(Q(R_{{\rm Ric}_0}))
  =     -\frac{2}{n-2} ({\rm Ric}_0^2)_0 + \frac{|{\rm Ric}_0|^2}{n} {\rm id}.
\]
Therefore, we have 
\begin{equation}\label{Tri-R-Full-Cal-Aux-7}
\begin{split}
\langle Q(R_{{\rm Ric}_0}), R_I \rangle&= \langle  \frac{|{\rm Ric}_0|^2}{n(n-1)} I,   \frac{\bar{\lambda}}{n-1} I\rangle
=\frac{1}{2(n-1)}\bar{\lambda}|{\rm Ric}_0|^2,
\end{split}
\end{equation}
and
\begin{equation}\label{Tri-R-Full-Cal-Aux-9}
\begin{split}
\langle Q(R_{{\rm Ric}_0}), R_{{\rm Ric}_0} \rangle&=\langle  \frac{2}{n-2}{\rm Ric}_0(R_{{\rm Ric}_0}^2+R_{{\rm Ric}_0}^\sharp) \wedge {\rm I},  \frac{2}{n-2}{\rm Ric}_0 \wedge {\rm I}\rangle\\
&=\Big(\frac{2}{n-2}\Big)^2\frac{n-2}{4}\langle {\rm Ric}_0(R_{{\rm Ric}_0}^2+R_{{\rm Ric}_0}^\sharp), {\rm Ric}_0\rangle\\
&=-\frac{2}{(n-2)^2}\langle  ({\rm Ric}_0^2)_0, {\rm Ric}_0\rangle.
\end{split}
\end{equation}
Simarly, by simply applying \eqref{BW-Q-W-Identity-2} in Lemma \ref{BW-Riccitype-Cur-1}, we obtain
\begin{equation}\label{Tri-R-Full-Cal-Aux-11}
\begin{split}
\langle Q(R_{{\rm Ric}_0}), W \rangle
= \frac{1}{n-2}\langle{\rm Ric}_0\wedge {\rm Ric}_0,W\rangle.
\end{split}
\end{equation}
Finally, with Lemma \ref{Q-R-Kahler-Decomposition-1}, we conclude that 
\begin{equation}\label{Tri-R-Full-Cal-Aux-13}
\begin{split}
\langle Q(R), W \rangle&=\frac{1}{2}{\rm tri}(W,R,R)\\
&=\langle Q(W,R_I),R\rangle +\langle Q(W,R_{{\rm Ric}_0}),R\rangle+\langle Q(W,W),R\rangle\\
&=\langle Q(W,R_{{\rm Ric}_0}),R_{{\rm Ric}_0}\rangle+\langle Q(W,W),W\rangle\\
&=\langle Q(R_{{\rm Ric}_0}), W \rangle+\langle Q(W), W \rangle.
\end{split}
\end{equation}
Taking the summation of \eqref{Tri-R-Full-Cal-Aux-3}, \eqref{Tri-R-Full-Cal-Aux-5} and \eqref{Tri-R-Full-Cal-Aux-13}, plugging \eqref{Tri-R-Full-Cal-Aux-7}  \eqref{Tri-R-Full-Cal-Aux-9} and  \eqref{Tri-R-Full-Cal-Aux-11}
we conclude that 

\begin{equation}
\begin{split}
\langle Q(R), R \rangle&=
\langle Q(R),R_I\rangle+\langle Q(R),R_{{\rm Ric}_0}\rangle+\langle Q(R),W\rangle\\
&=\langle Q(R),R_I\rangle+2\langle Q(R_{{\rm Ric}_0}), R_I \rangle +3\langle Q(R_{{\rm Ric}_0}), W \rangle\\
&\quad+\langle Q(R_{{\rm Ric}_0}), R_{{\rm Ric}_0} \rangle+\langle Q(W), W \rangle\\
&=\frac{\bar{\lambda}}{2(n-1)}\big(|{\rm Ric}_0|^2+n\bar{\lambda}^2\big)+\frac{1}{n-1}\bar{\lambda}|{\rm Ric}_0|^2+\frac{3}{n-2}\langle{\rm Ric}_0\wedge {\rm Ric}_0,W\rangle\\
&\quad -\frac{2}{(n-2)^2}\langle  ({\rm Ric}_0^2)_0, {\rm Ric}_0\rangle+\langle Q(W), W\rangle.
\end{split}
\end{equation}
The proof is finished.
\end{proof}
\begin{lemma} For $R\in S_B^2(\mathfrak{so}(n))$, we have  
\begin{equation}\label{R-Ric-2-Inn-Pro-For-1}
\begin{split}
\big\langle {\rm Ric}(R), {\rm Ric}(Q(R))\big\rangle&=2\langle R,{\rm Ric}\wedge {\rm Ric}\rangle\\
&=n\bar{\lambda}^3+\frac{2n-3}{n-1}\bar{\lambda}|{\rm Ric}_0|^2-
    \frac{2}{n-2}\langle{\rm Ric}_0, ({\rm Ric}_0^2)_0 \rangle\\
&\quad+2\langle W,  {\rm Ric}_0\wedge {\rm Ric}_0\rangle.
\end{split}
\end{equation}
\end{lemma}
\begin{proof} The first equality follows from  \eqref{Ricci-QR-For-1} and 
 \begin{align*}
 \langle R,{\rm Ric}\wedge {\rm Ric}\rangle&=\frac{1}{2} \sum_{i,j}\langle R(e_i\wedge e_j),{\rm Ric}\wedge {\rm Ric}(e_i\wedge e_j)\rangle\\
 &=\frac{1}{4} \sum_{i,j,k,l}  \big(R_{ijkl}R_{ik}R_{jl}- R_{ijkl}R_{il}R_{jk}\big)\\
 &= \frac{1}{2}\langle {\rm Ric}, \mathring{R}({\rm Ric})\big\rangle.
 \end{align*}
Furthermore,  with  \eqref{Curvature-Decomposition-1},  we calculate 
\[\begin{split}
\langle R,  {\rm Ric}_0\wedge {\rm Ric}_0\rangle&= \langle \frac{\bar{\lambda}}{n-1}{\rm id} \wedge {\rm id}
    +\frac{2}{n-2}{\rm Ric}_0 \wedge {\rm id}+W, {\rm Ric}_0\wedge {\rm Ric}_0\rangle \\
        &=-\frac{\bar{\lambda}}{2(n-1)}|{\rm Ric}_0|^2-
    \frac{1}{n-2}\langle{\rm Ric}_0, ({\rm Ric}_0^2)_0 \rangle+\langle W, {\rm Ric}_0\wedge {\rm Ric}_0\rangle.
\end{split}  \]
We proceed as above by using \eqref{Curvature-Decomposition-1}    once again leads to
\[\begin{split}
\langle R,   {\rm Ric}\wedge  {\rm Ric}\rangle&=\langle R,  (\bar{\lambda}{\rm id}+{\rm Ric}_0)\wedge (\bar{\lambda}{\rm id}+{\rm Ric}_0)\rangle\\
&=\bar{\lambda}^2\langle R_I, I\rangle + 2\bar{\lambda}\langle R_{{\rm Ric}_0},   {\rm Ric}_0\wedge {\rm id}\rangle+\langle R, {\rm Ric}_0 \wedge {\rm Ric}_0\rangle\\
&=\frac{n}{2}\bar{\lambda}^3+\frac{2n-3}{2(n-1)}\bar{\lambda}|{\rm Ric}_0|^2-
    \frac{1}{n-2}\langle{\rm Ric}_0, ({\rm Ric}_0^2)_0 \rangle+\langle W,  {\rm Ric}_0\wedge {\rm Ric}_0\rangle.
\end{split}     \]
The proof is finished.
\end{proof}

\subsection{Some algebraic  inequalities for curvature operators}
Now we  collect some algebraic  inequalities for curvature operators which will be used in our proof. 

First, we need the following estimates for symmetric tensors. 
\begin{lemma}{\cite[Lemma 2.4]{MR0806701}}\label{Huisken-Basic-Inequality}
Let $T=(T_{ij})_{~1\leq i, j\leq m}$ be a symmetric trace-free operator with eigenvalues $\lambda_1,\cdots, \lambda_m$ such that 
\[\sum_{i=1}^m\lambda_i=0,\quad |T|^2=\sum_{i=1}^m\lambda_i^2.\]
Then we have 
\begin{equation}\label{Trace-Free-Eigenvalue-Estimate-1}
\lambda_i^2\leq \frac{m-1}{m}\|T\|^2, ~~1\leq i \leq m,
\end{equation}
and
\begin{equation}\label{Huisken-3rd-Inequality-1}
|{\rm tr}(T^3)|=\Big|\sum_{i=1}^m\lambda_i^3\Big|\leq \frac{m-2}{\sqrt{m(m-1)}}\|T\|^3.
\end{equation}
The equality in \eqref{Huisken-3rd-Inequality-1} holds only  if either 
\[\lambda_{1}=\cdots=\lambda_{m}=0,\]
or (up to permutation and sign)
\[
\lambda_{1}=\cdots=\lambda_{m-1}=-\frac{1}{\sqrt{m(m-1)}}\|T\|,\quad \lambda_m= \sqrt{\frac{m-1}{m}}\|T\|.
\]
\end{lemma}
Apply the  lemma  \ref{Huisken-Basic-Inequality} to ${\rm Ric}_0$, we immediately obtain from  \eqref{Trace-Free-Eigenvalue-Estimate-1}  that 
\begin{corollary}\label{Ricci-Positive-Cond-1}
  If $a>\frac{n}{4}-\frac{1}{2}$, then we have 
  \[\Omega(a)\subset \Big\{R\in S_B^2(\mathfrak{so}(n))\,\big|\,{\rm Ric}(R)>0\Big\}.\]
\end{corollary}
\begin{proof} If $R\in \Omega(a)$, see \eqref{B-W-Cone-Def-1},  we have 
\[a \|R_{{\rm Ric}_0}\|^2\leq \frac{n-4a}{4}\|R_I\|^2,\]
it is rewritten equivalently by    \eqref{curvature-norm-decomposition-3}
\eqref{curvature-norm-decomposition-5} as  
\[a \frac{1}{n-2}|{\rm Ric}_0|^2\leq \frac{n-4a}{4}\frac{n}{2(n-1)}\bar{\lambda}^2.\]
Hence, we conclude  from   \eqref{Trace-Free-Eigenvalue-Estimate-1}  that 
\[\lambda_i^2({\rm Ric}_0)\leq \frac{n-1}{n}|{\rm Ric}_0|^2\leq \frac{n-4a}{4a}\frac{n-2}{2}\bar{\lambda}^2,~\forall i.\]
Therefore,  if $\frac{n-4a}{4a}\frac{n-2}{2}<1$, i.e., $a>\frac{n}{4}-\frac{1}{2}$, we have  ${\rm Ric}(R)>0$. 
\end{proof}

Apply the  lemma  \ref{Huisken-Basic-Inequality} to ${\rm Ric}_0$ again, we immediately obtain from \eqref{Huisken-3rd-Inequality-1}  that 
\begin{proposition}\cite[in proof of  Theorem 3.3]{MR0806701}
\begin{equation}\label{Huisken-Trace-Free-Inequality-1}
\begin{split}
\big|\langle  ({\rm Ric}_0^2)_0, {\rm Ric}_0\rangle\big| \leq \frac{n-2}{\sqrt{n(n-1)}}|{\rm Ric}_0|^3.
\end{split}
\end{equation}
\end{proposition}
The tensor ${\rm Ric}_0\wedge {\rm Ric}_0$ has the same symmetries as the algebraic Riemannian tensor,   therefore by dissociating  the `Weyl' parts  from the   scalar  and  traceless Ricci part, we have 
\begin{lemma}{\cite[Lemma 3.4]{MR0806701}}  
\begin{equation}\label{Huisken-Weyl-Pinching-1}
\big|\langle{\rm Ric}_0\wedge {\rm Ric}_0,W\rangle\big|\leq \sqrt{\frac{n-2}{2(n-1)}}\|W\||{\rm Ric}_0|^2.
\end{equation}
\end{lemma}

The following estimate  was  derived  by Huisken in 
his proof of the Theorem 3.3 in \cite{MR0806701}. 
\begin{lemma}For $n\geq 4$, we have  \cite[in proof of  Theorem 3.3]{MR0806701}
\begin{equation}\label{Huisken-Trace-W-Inequality-1}
\begin{split}
|\langle Q(W), {W}\rangle|&\leq  \sqrt{\frac{(n^2-1)(n-2)}{n}}\|W\|^3.
\end{split}
\end{equation}
There is an improvement in dimension 4 (see [Lemma 3.5]{\cite{MR0806701}}). With the Lie algebraic decomposition $\mathfrak{so}(4)=\mathfrak{so}(3)\oplus \mathfrak{so}(3)$,  the Weyl tensor may be split into two parts,  
and thus: 
\begin{equation}\label{Huisken-Trace-W-Inequality-n4}
\begin{split}
|\langle Q(W), {W}\rangle|&\leq  \frac{\sqrt{6}}{2}\|W\|^3.
\end{split}
\end{equation}
\end{lemma}
\begin{proof} In  case of $n=4$,   the estimate  \eqref{Huisken-Trace-W-Inequality-n4}  follows from  [Lemma 3.5]{\cite{MR0806701}}. 
In case of $n\geq 5$,   since the estimate \eqref{Huisken-Trace-W-Inequality-1} was not written down as an independent lemma in \cite{MR0806701},     we include the calculations here for readers' convenience.
 
Using an idea of Tachibana \cite{MR0365415}, for each fixed $i,j,k,l$, we define a local skew-symmetric $2$-tensor field $\xi^{ijkl}=(\xi^{pq}_{ijkl})$ by 
\[\begin{split}\xi_{ijkl}&=e_i\wedge W(e_k,e_l, e_j)+W(e_k,e_l, e_i)\wedge e_j  +e_k\wedge W(e_i,e_j, e_l)+ W(e_i,e_j, e_k)\wedge e_l.
\end{split}\]
A straightforward computation then gives 
\begin{align}
\langle W(\xi_{ijkl}),\xi_{ijkl} \rangle
&= 4 W_{pjkl}W_{qjkl}W_{piqi} -4W_{ipkl}W_{qjkl}W_{iqp j}-8W_{ijpl}W_{qjkl}W_{iqpk}.\label{Q-W-Huisken-Estimate-Aux-3}
\end{align}
Using then the Bainchi identities, 
\[\begin{split}
 W_{ipkl}W_{qjkl}W_{iqpj} 
&= W_{ipkl}W_{qjkl}W_{ipqj}-W_{ipkl}W_{jqkl}W_{ijpq}, 
\end{split}\]
 thus 
\[
 W_{ipkl}W_{qjkl}W_{iqpj}=\frac{1}{2} W_{ipkl}W_{klqj}W_{qjip},\]
we obtain from \eqref{Q-W-Huisken-Estimate-Aux-3} that 
\begin{equation}\label{Q-W-Huisken-Estimate-Aux-5}
\begin{split}
\langle W(\xi_{ijkl}),\xi_{ijkl} \rangle
&=-2( W_{ipkl}W_{klqj}W_{qjip}+4W_{ijpl}W_{qjkl}W_{iqpk}).
\end{split}
\end{equation}
On the other hand,  with \eqref{Cur-Ope-Squ-Loc-For-1} and  \eqref{Cur-Ope-Sharp-Loc-For-1}, we have 
\begin{equation}\label{Q-W-Huisken-Estimate-Aux-7}
\begin{split}
\langle Q(W), W\rangle&=\langle W^2+W^\sharp, W\rangle\\
&=\frac{1}{8}\sum_{i,j,k,l}\Big(\sum_{p,q=1}^nW_{ijpq}W_{klpq}W_{ijkl}+4\sum_{p,q=1}^nW_{ipkq}W_{jplq}W_{ijkl}\Big).
\end{split}
\end{equation}
By comparing \eqref{Q-W-Huisken-Estimate-Aux-5} with  \eqref{Q-W-Huisken-Estimate-Aux-7}, we conclude that 
\begin{equation}\label{Q-W-Estimate-Aux-3}\sum_{p,q,r,s}\langle  W(\xi_{pqrs}), \xi_{pqrs}\rangle=-16\langle Q(W), W\rangle. \end{equation}
Moreover, a direct calculation shows that  
\begin{align*}
\langle  \xi_{ijkl},\xi_{ijkl} \rangle&=\langle  W_{pjkl}e_p\wedge e_i +W_{ipkl}e_p\wedge e_j +W_{ijpl}e_p\wedge e_k +W_{ijkp}e_p\wedge e_l,\\
&\qquad W_{qjkl}e_q\wedge e_i +W_{iqkl}e_q\wedge e_j +W_{ijql}e_q\wedge e_k +W_{ijkq}e_q\wedge e_l \rangle\\
&= 4(n-1)  W_{ijkl}W_{ijkl}-4 W_{ipkl}W_{ipkl}-8W_{ijpl}W_{iljp}.
\end{align*}
Note that by Bainchi identities again, 
\[
W_{ijpl}W_{iljp}=-W_{ijpl}(W_{ljip}+W_{jilp}),\]
we conclude that 
\[W_{ijpl}W_{iljp} =-\frac{1}{2}W_{ijpl}W_{ijpl}.
\]
Therefore
\begin{equation}\label{Q-W-Estimate-Aux-5}
\sum_{p,q,r,s}\langle \xi^{pqrs}, \xi^{pqrs}\rangle=4(n-1)|W|^2=16(n-1)\|W\|^2.
\end{equation}

 Since $W$ is a symmetric trace-free operator on the space of skew-symmetric $2$-tensor field, then we have the estimate from  \eqref{Trace-Free-Eigenvalue-Estimate-1} in Lemma \ref{Huisken-Basic-Inequality},
\begin{equation}\label{Q-W-Estimate-Aux-7}
\Big|\frac{\langle W(\xi^{pqrs}), \xi^{pqrs}\rangle}{\langle \xi^{pqrs}, \xi^{pqrs}\rangle}\Big|\leq \sqrt{\frac{N-1}{N}}\|W\|,\quad N=\frac{n(n-1)}{2}.
\end{equation}
Therefore, we conclude  from  \eqref{Q-W-Estimate-Aux-3},  \eqref{Q-W-Estimate-Aux-5} and \eqref{Q-W-Estimate-Aux-7}  that 
\begin{equation*}
\begin{split}
|\langle Q(W), W\rangle |&\leq \frac{1}{16}\sum_{p,q,r,s}|\langle W(\xi^{pqrs}), \xi^{pqrs}\rangle|\\
&\leq  \frac{1}{16}\sqrt{\frac{N-1}{N}}\|W\|\sum_{p,q,r,s}\langle \xi^{pqrs}, \xi^{pqrs}\rangle\\
&=\sqrt{\frac{(n^2-1)(n-2)}{n}}\|W\|^3.
\end{split}
\end{equation*}
The conclusion of the lemma is proved.
\end{proof}

\section{Ricci flow invariant curvature conditions and sphere theorem}\label{Proof-nain-theorem}
Now, we begin to  prove Theorems   \ref{Main-Theorem-1} and \ref{Main-Convergence-Thm}. 

\subsection{Proof of  Theorem \ref{Main-Theorem-1}}
\begin{proof}[Proof of  Theorem \ref{Main-Theorem-1}]
On the set ${\rm scal}_+=\{R\in S_B^2(\mathfrak{so}(n))\, |{\rm scal}(R)>0\}$,  we can define the $O(n)$-invariant function
\[F:{\rm scal}_+ \rightarrow \mathbb{R},~R\mapsto \frac{\|R\|^2}{|{\rm Ric}(R)|^2}.\]
Since  ${\rm Ric}(R)\neq 0$, for all $R\in {\rm scal}_+$, then the function $F$ is well-defined. Restricted to the cone $\Omega(0) =\Big\{R\in S_B^2(\mathfrak{so}(n))\,\Big| \frac{n-2}{4} \|R\|^2\leq \frac{1}{4}|{\rm Ric}|^2~\text{and}~{\rm scal}(R)>0\Big\}$, 
$F$ is bounded. More specifically, 
\[\frac{1}{2(n-1)}=F(I)\leq F(R)\leq\frac{1}{n-2}, ~~ R\in \Omega.\] Moreover, $F$ is smooth and for each $R\in {\rm scal}_+$ and its differential is given by 
\begin{equation}\label{Diff-F-For-1}
dF_R(S) = \frac{2}{|{\rm Ric}(R)|^2}\Big(\big\langle R, S\big\rangle -F(R)\big\langle {\rm Ric}(R), {\rm Ric}(S)\big\rangle \Big),\end{equation}
for all $S\in T_R{\rm scal}_+\cong S_B^2(\mathfrak{so}(n)).$

If we take $S=Q(R)=R^2+R^\sharp$ in \eqref{Diff-F-For-1},   we obtain that 
\begin{equation}\label{Diff-F-For-Q-1}
dF_R(Q(R)) = \frac{2}{|{\rm Ric}(R)|^2}\Big(\big\langle R, Q(R)\big\rangle -F(R)\big\langle {\rm Ric}(R), {\rm Ric}(Q(R))\big\rangle \Big).
\end{equation}
By our construction, note that  $\partial\Omega(a)$ is the level set $F^{-1}(\frac{1}{n-2+4a})$. With Proposition \ref{Hamilton'-maximum-principle} and \ref{RF-Inv-Tangent-Cone-1},  to prove   the main Theorem \ref{Main-Theorem-1},
 we only need to verify that
\begin{equation}\label{Diff-F-For-Q-neg-1}
dF_R(Q(R))\leq 0,\quad  \forall R\in \partial\Omega(a).
\end{equation}
Hence, the main task is to   look  for all the possible  $a  \in [0, n/4]$  such that (see \eqref{Diff-F-For-Q-1}), 
\begin{equation}\label{Diff-F-For-Q-neg-3}
\big\langle  R,Q(R)\big\rangle -F(R)\big\langle {\rm Ric}(R), {\rm Ric}(Q(R))\big\rangle \leq 0, \quad  \forall R\in \partial\Omega(a).
\end{equation}
Now, by \eqref{Tri-R-Full-Huisken-1} and \eqref{R-Ric-2-Inn-Pro-For-1}, we have
\begin{align}
\big\langle& R,Q(R)\big\rangle -F(R)\big\langle {\rm Ric}(R), {\rm Ric}(Q(R))\big\rangle\notag\\
    &= \frac{n\bar{\lambda}^3}{2(n-1)} +\frac{3}{2(n-1)}\bar{\lambda}|{\rm Ric}_0|^2+\frac{3}{n-2}\langle{\rm Ric}_0\wedge {\rm Ric}_0,W\rangle\notag\\
&\quad -\frac{2}{(n-2)^2}\langle  ({\rm Ric}_0^2)_0, {\rm Ric}_0\rangle+ \langle Q(W),W \rangle
\notag\\
    &\quad    -F(R)\Big(n\bar{\lambda}^3+\frac{2n-3}{n-1}\bar{\lambda}|{\rm Ric}_0|^2-
    \frac{2}{n-2}\langle{\rm Ric}_0, ({\rm Ric}_0^2)_0 \rangle+2\langle W,  {\rm Ric}_0\wedge {\rm Ric}_0\rangle
\Big)\notag\\
&=\Big(\frac{1}{2(n-1)}-F(R)\Big)n\bar{\lambda}^3+\Big(\frac{3}{2(2n-3)}-F(R)\Big)\frac{2n-3}{n-1}\bar{\lambda}|{\rm Ric}_0|^2\notag\\
&\quad  -\Big(\frac{1}{n-2}-F(R)\Big)\frac{2}{n-2}\langle{\rm Ric}_0, ({\rm Ric}_0^2)_0 \rangle \label{Tri-R-Ric-Wilking-3}\\
&\quad +2\Big(\frac{3}{2(n-2)}-F(R)\Big)\langle{\rm Ric}_0\wedge {\rm Ric}_0,W\rangle+\langle Q(W),W \rangle. \notag
\end{align}
Fix $a>0$,  at $ R\in \partial\Omega(a)=F^{-1}(\frac{1}{n-2+4a})$, a direct calculation shows that 
\begin{align*}
\frac{1}{2(n-1)}-F(R)&=\frac{1}{2(n-1)}-\frac{1}{n-2+4a}=\frac{4a-n}{2(n-1)(n-2+4a)}\leq 0,\\
\frac{3}{2(2n-3)}-F(R)&=\frac{3}{2(2n-3)}-\frac{1}{n-2+4a}=\frac{12a-n}{2(2n-3)(n-2+4a)},\\
\frac{1}{n-2}-F(R)&=\frac{1}{n-2}-\frac{1}{n-2+4a}=\frac{4a}{(n-2)(n-2+4a)}\geq 0,\\
\frac{3}{2(n-2)}-F(R)&=\frac{3}{2(n-2)}-\frac{1}{n-2+4a}=\frac{12a+n-2}{2(n-2)(n-2+4a)}\geq 0.
\end{align*}
In particular, the coefficient of the first term in  \eqref{Tri-R-Ric-Wilking-3} have favorite  sign, and we will see this  negative term will dominate all the other terms when $a$ is close to $n/4$. Moreover,  
the coefficients of the third and fourth terms in  \eqref{Tri-R-Ric-Wilking-3} have fixed sign, 
therefore,   from \eqref{Tri-R-Ric-Wilking-3},  we have 
\begin{equation}\label{Tri-R-Ric-Wilking-5}
\begin{split}
\big\langle& R,Q(R)\big\rangle -F(R)\big\langle {\rm Ric}(R), {\rm Ric}(Q(R))\big\rangle\\
&\leq\frac{4a-n}{2(n-1)(n-2+4a)}n\bar{\lambda}^3+\frac{12a-n}{2(n-1)(n-2+4a)}\bar{\lambda}|{\rm Ric}_0|^2\\
&\quad +\frac{8a}{(n-2)^2(n-2+4a)}\big|\langle{\rm Ric}_0, ({\rm Ric}_0^2)_0 \rangle\big|\\
&\quad+\frac{12a+n-2}{(n-2)(n-2+4a)}\big|\langle{\rm Ric}_0\wedge {\rm Ric}_0,W\rangle\big| +\big|\langle Q(W),W \rangle\big|. 
\end{split}
\end{equation}
Next, plug  the  estimates  \eqref{Huisken-Trace-Free-Inequality-1}, \eqref{Huisken-Weyl-Pinching-1}  and \eqref{Huisken-Trace-W-Inequality-1} into \eqref{Tri-R-Ric-Wilking-5}, we  conclude that 
\begin{equation}\label{Tri-R-Ric-Wilking-7}
\begin{split}
\big\langle R,& Q(R)\big\rangle -F(R)\big\langle {\rm Ric}(R), \mathring{R}({\rm Ric})\big\rangle\\
&\leq  \frac{4a-n}{2(n-1)(n-2+4a)}n\bar{\lambda}^3+\frac{12a-n}{2(n-1)(n-2+4a)}\bar{\lambda}|{\rm Ric}_0|^2\\
&\quad  +\frac{8a}{(n-2)(n-2+4a)}\frac{1}{\sqrt{n(n-1)}}|{\rm Ric}_0|^3\\
&\quad +\frac{12a+n-2}{(n-2)(n-2+4a)}\frac{\sqrt{(n-2)}}{\sqrt{2(n-1)}}\|W\||{\rm Ric}_0|^2\\
&\quad +\sqrt{\frac{(n^2-1)(n-2)}{n}}\|W\|^3.
\end{split}
\end{equation}
At $R\in \partial\Omega(a)=F^{-1}(\frac{1}{n-2+4a})$, from \eqref{B-W-Cone-Def-1},  we have
\[ a \|R_{{\rm Ric}_0}\|^2+\Big(\frac{n-2}{4}+a\Big)\|W\|^2=\frac{n-4a}{4}\|R_I\|^2,
 \]
 or equivalently,    
 \[\frac{a}{n-2} |{\rm Ric}_0|^2+\Big(\frac{n-2}{4}+a\Big)\|W\|^2= \frac{n-4a}{4}\frac{n}{2(n-1)}\bar{\lambda}^2.\]
We may set
\[|{\rm Ric}_0|^2=\frac{n-2}{a}\frac{n-4a}{4}\frac{n}{2(n-1)}\bar{\lambda}^2\cos^2\theta=\frac{n(n-2)(n-4a)}{8a(n-1)}\bar{\lambda}^2\cos^2\theta,\]
and
\[\|W\|^2=\frac{4}{n-2+4a}\frac{n-4a}{4}\frac{n}{2(n-1)}\bar{\lambda}^2\sin^2\theta=\frac{n(n-4a)}{2(n-1)(n-2+4a)}\bar{\lambda}^2\sin^2\theta,\]
then the inequality  \eqref{Tri-R-Ric-Wilking-7} can be rewritten as 
\begin{align}
\big\langle R,& Q(R)\big\rangle -F(R)\big\langle {\rm Ric}(R), \mathring{R}({\rm Ric})\big\rangle\notag\\
&\leq  
 \frac{n(4a-n)}{2(n-1)(n-2+4a)}\bar{\lambda}^3+\frac{n(n-2)}{16(n-1)^2}\big(8+\frac{4a-n}{a}\big)\frac{n-4a}{n-2+4a}\cos^2\theta\bar{\lambda}^3\notag\\
&+\frac{n\sqrt{n-2}}{2\sqrt{2}(n-1)^2}\frac{n-4a}{(n-2+4a)}\frac{\sqrt{n-4a}}{\sqrt{a}}|\cos\theta|^3\bar{\lambda}^3\notag\\
&+\frac{12a+n-2}{16a}\frac{n\sqrt{n(n-2)}}{(n-1)^2}\frac{n-4a}{n-2+4a}\frac{\sqrt{n-4a}}{\sqrt{n-2+4a}}\cos^2\theta|\sin\theta|\bar{\lambda}^3\notag\\
&+\frac{n\sqrt{(n+1)(n-2)}}{2\sqrt{2}(n-1)}\frac{n-4a}{n-2+4a}\frac{\sqrt{n-4a}}{\sqrt{n-2+4a}}|\sin\theta|^3\bar{\lambda}^3\notag\\
&=  \frac{(4a-n)\bar{\lambda}^3}{n-2+4a}\frac{n}{2(n-1)}\Big(\sin^2\theta+\frac{1}{n-1}\cos^2\theta+\frac{n-2}{8(n-1)}\frac{n-4a}{a}\cos^2\theta\notag\\
&-\frac{\sqrt{n-2}}{\sqrt{2}(n-1)}\frac{\sqrt{n-4a}}{\sqrt{a}}|\cos\theta|^3-\frac{\sqrt{n(n-2)}}{n-1}\frac{12a+n-2}{8a}\frac{\sqrt{n-4a}}{\sqrt{n-2+4a}}\cos^2\theta|\sin\theta|\label{Diff-F-For-Q-nonnegative-4}\\
&-\frac{\sqrt{(n+1)(n-2)}}{\sqrt{2}}\frac{\sqrt{n-4a}}{\sqrt{n-2+4a}}|\sin\theta|^3\Big).\notag
 \end{align}

Using then the H\"older inequality, 
  \[|\sin\theta| \leq \epsilon +\frac{1}{4\epsilon} \sin^2\theta=\epsilon +\frac{1}{4\epsilon} (1-\cos^2\theta),\]
 \[|\cos\theta| \leq \eta +\frac{1}{4\eta} \cos^2\theta=\eta +\frac{1}{4\eta} (1-\sin^2\theta),\]
we obtain from \eqref{Diff-F-For-Q-nonnegative-4} that 
 \begin{align}
\big\langle R,& Q(R)\big\rangle -F(R)\big\langle {\rm Ric}(R), \mathring{R}({\rm Ric})\big\rangle\notag\\
&\leq   \frac{(4a-n)\bar{\lambda}^3}{n-2+4a}\frac{n}{2(n-1)}\Big(\sin^2\theta+\frac{1}{n-1}\cos^2\theta+\frac{n-2}{8(n-1)}\frac{n-4a}{a}\cos^2\theta\notag\\
&-\frac{\sqrt{n-2}}{\sqrt{2}(n-1)}\frac{\sqrt{n-4a}}{\sqrt{a}}\cos^2\theta\big(\eta +\frac{1}{4\eta} (1-\sin^2\theta)\big)\notag\\
&-\frac{\sqrt{n(n-2)}}{n-1}\frac{12a+n-2}{8a}\frac{\sqrt{n-4a}}{\sqrt{n-2+4a}}\cos^2\theta\big(\epsilon_1 +\frac{1}{4\epsilon_1} \sin^2\theta\big)\notag\\
&-\frac{\sqrt{(n+1)(n-2)}}{\sqrt{2}}\frac{\sqrt{n-4a}}{\sqrt{n-2+4a}}\sin^2\theta\big(\epsilon_2 +\frac{1}{4\epsilon_2} (1-\cos^2\theta)\big)\Big)\notag\\
&=  \frac{(4a-n)\bar{\lambda}^3}{n-2+4a}\frac{n}{2(n-1)}\Big(\sin^2\theta\Big[1-\frac{\sqrt{(n+1)(n-2)}}{\sqrt{2}}\frac{\sqrt{n-4a}}{\sqrt{n-2+4a}}(\epsilon_2 +\frac{1}{4\epsilon_2})\Big]\notag\\
&+\cos^2\theta \frac{1}{n-1}\Big[1+\frac{n-2}{8}\frac{n-4a}{a}-\frac{\sqrt{n-2}}{\sqrt{2}}\frac{\sqrt{n-4a}}{\sqrt{a}}(\eta +\frac{1}{4\eta})\notag\\
&-\sqrt{n(n-2)}\frac{12a+n-2}{8a}\frac{\sqrt{n-4a}}{\sqrt{n-2+4a}} \epsilon_1\Big]\label{Diff-F-For-Q-nonnegative-5},\\
&+\sin^2\theta\cos^2\theta\frac{\sqrt{n-2}}{4\sqrt{2}(n-1)}\frac{\sqrt{n-4a}}{\sqrt{n-2+4a}}\theta\Big[\frac{\sqrt{n-2+4a}}{\sqrt{a}}\frac{1}{\eta}\notag\\
&-\sqrt{2n}\frac{12a+n-2}{8a}\frac{1}{\epsilon_1}+(n-1)\sqrt{n+1}\frac{1}{\epsilon_2}\Big] \Big).\notag
 \end{align}
 
To have  \eqref{Diff-F-For-Q-neg-3} holds, by \eqref{Diff-F-For-Q-nonnegative-5},  we  are now  looking for the smallest $a\in (0, n/4]$  with  appropriate choice of $\eta$, $\epsilon_1$, $\epsilon_2$ to guarantee   the following three inequalities hold simultaneously  
 \begin{equation}\label{Aux-a-n-sin-1}
 1-\frac{\sqrt{(n+1)(n-2)}}{\sqrt{2}}\frac{\sqrt{n-4a}}{\sqrt{n-2+4a}}(\epsilon_2 +\frac{1}{4\epsilon_2})\geq 0,
 \end{equation}
 \begin{equation}\label{Aux-a-n-cos-2}
 \begin{split}1&+\frac{n-2}{8}\frac{n-4a}{a}-\frac{\sqrt{n-2}}{\sqrt{2}}\frac{\sqrt{n-4a}}{\sqrt{a}}(\eta +\frac{1}{4\eta})\\
 & -\sqrt{n(n-2)}\frac{12a+n-2}{8a}\frac{\sqrt{n-4a}}{\sqrt{n-2+4a}} \epsilon_1\geq 0,\end{split}
  \end{equation}
 \begin{equation}\label{Aux-a-n-sincos-3}
 \frac{\sqrt{n-2+4a}}{\sqrt{a}}\frac{1}{\eta}-\sqrt{2n}\frac{12a+n-2}{8a}\frac{1}{\epsilon_1}+(n-1)\sqrt{n+1}\frac{1}{\epsilon_2}\geq 0.
  \end{equation}
 First, the inequality \eqref{Aux-a-n-sin-1} forces 
  \begin{equation}\label{Aux-a-n-sin-1-1}
 \epsilon_2 +\frac{1}{4\epsilon_2}\leq  \frac{\sqrt{2}\sqrt{n-2+4a}}{\sqrt{(n+1)(n-2)}\sqrt{n-4a}};
 \end{equation}
If
 \begin{equation}\label{Aux-a-n-0-0}
 \frac{\sqrt{2}\sqrt{n-2+4a}}{\sqrt{(n+1)(n-2)}\sqrt{n-4a}}\geq 1, 
  \end{equation}
 i.e.,
 \begin{equation}\label{Aux-a-n-0-1}
 a\geq \frac{n^2-4}{4n}=\frac{n}{4}-\frac{1}{n},
 \end{equation}
it is possible to choose 
  \begin{equation}\label{Aux-a-n-sin-eps-2}
  \epsilon_2= \frac{\sqrt{(n+1)(n-2)}\sqrt{n-4a}}{2\sqrt{2}\sqrt{n-2+4a}}\leq \frac{1}{2},
  \end{equation} 
such that \eqref{Aux-a-n-sin-1-1} holds, and  then  \eqref{Aux-a-n-sin-1}  follows immediately from \eqref{Aux-a-n-sin-1-1}.

Next, we choose \begin{equation}\label{eta--def-aux-1}
\eta=\frac{\sqrt{n-2}\sqrt{n-4a}}{2\sqrt{2}\sqrt{a}},
\end{equation}
then  \eqref{Aux-a-n-sincos-3} can be rewritten as

 \begin{equation}\label{Aux-a-n-sincos-3-1}
 \begin{split}
& \frac{\sqrt{n-2+4a}}{\sqrt{a}}\frac{2\sqrt{2}\sqrt{a}}{\sqrt{n-2}\sqrt{n-4a}}-\sqrt{2n}\frac{12a+n-2}{8a}\frac{1}{\epsilon_1}\\
 &\quad+(n-1)\sqrt{n+1}\frac{2\sqrt{2}\sqrt{n-2+4a}}{\sqrt{(n+1)(n-2)}\sqrt{n-4a}}\geq 0,
 \end{split}
  \end{equation} 
i.e., 
\[
-\frac{12a+n-2}{8a}\frac{1}{\epsilon_1}+ 2\sqrt{n}\frac{\sqrt{n-2+4a}}{\sqrt{n-2}\sqrt{n-4a}}\geq 0.
\]
Therefore,  it is possible to choose 
  \begin{equation}\label{Aux-a-n-sin-eps-1}
  \epsilon_1= \frac{12a+n-2}{16\sqrt{n}a}
\frac{\sqrt{n-2}\sqrt{n-4a}}{\sqrt{n-2+4a}},
\end{equation}
such that \eqref{Aux-a-n-sincos-3-1}, as well as  \eqref{Aux-a-n-sincos-3} hold. 

Choose  $\epsilon_1, \epsilon_2$  and $\eta$ as in \eqref{Aux-a-n-sin-eps-1}, \eqref{Aux-a-n-sin-eps-2} and \eqref{eta--def-aux-1} respectively,  we can rewrite the left hand side of \eqref{Aux-a-n-cos-2} as following
\begin{equation}\label{Aux-a-n-cos-2-3}
 \begin{split}
 1&+\frac{(n-2)(n-4a)}{8a}-\frac{\sqrt{n-2}}{\sqrt{2}}\frac{\sqrt{n-4a}}{\sqrt{a}}(\frac{\sqrt{n-2}\sqrt{n-4a}}{2\sqrt{2}\sqrt{a}}\\
 &\quad+\frac{\sqrt{a}}{\sqrt{2}\sqrt{n-2}\sqrt{n-4a}})-\frac{(12a+n-2)^2}{128a^2}
\frac{(n-2)(n-4a)}{n-2+4a}\\
  &= \frac{1}{2}-\Big(\frac{1}{8}+\frac{(12a+n-2)^2}{128(n-2+4a)a}\Big)
\frac{(n-2)(n-4a)}{a}.
 \end{split}
  \end{equation}
Note that \eqref{Aux-a-n-cos-2-3} is increasing with respect to $a>0$ when $n\geq 4$. In fact,  we have 
\begin{equation}\label{Aux-Term-Mono-1}
\begin{split}\frac{d}{da}\frac{(12a+n-2)^2}{128(n-2+4a)a}
&=\frac{12a+n-2}{128(n-2+4a)^2a^2}(n-2)(4a-n+2),\\
\end{split}\end{equation}
and
\begin{equation}\begin{split}\frac{d}{da}&\Big[\frac{1}{2}-\Big(\frac{1}{8}+\frac{(12a+n-2)^2}{128(n-2+4a)a}\Big)\frac{(n-2)(n-4a)}{a}\Big]\\
&=-\frac{12a+n-2}{128(n-2+4a)^2a^2}(n-2)(4a-n+2)\frac{(n-2)(n-4a)}{a}\\
&\quad+\frac{16(n-2+4a)a+(12a+n-2)^2}{128(n-2+4a)a}\frac{n(n-2)}{a^2}\\
&=\frac{n-2}{64(n-2+4a)^2a^3}\Big[64(8n-3) a^3 +16(9n+2)(n-2)a^2\\
&\quad+4(n-2)^2(6n+1)a +n(n-2)^3)\Big]\\
&>0.
\end{split}\end{equation}
Therefore, if \eqref{Aux-a-n-0-1} holds, we conclude that 
 \begin{equation}\label{Aux-a-n-cos-2-5}
 \begin{split}
 \eqref{Aux-a-n-cos-2-3}
  &\geq \frac{1}{2}-\Big(\frac{1}{8}+\frac{(12a+n-2)^2}{128(n-2+4a)a}\Big)
\frac{(n-2)(n-4a)}{a}\Big|_{a=\frac{n}{4}-\frac{1}{n}}\\
&=\frac{(n+1)(n+2)(n-10)-2}{2(n+1)(n+2)^2}\\
&>0,
 \end{split}
\end{equation}
when $n\geq 11$.
Consequently, when $n\geq 11$, \eqref{Aux-a-n-cos-2-3} is positive and thus \eqref{Aux-a-n-cos-2} holds provide $n/4-\frac{1}{n}\leq a\leq n/4$.

 When $4\leq n\leq 10$, by the monotonicity of  \eqref{Aux-a-n-cos-2-3}, we may  estimate a lower bound of $a$ by solving $\eqref{Aux-a-n-cos-2-3}=0$. Alternatively, we may  give a rough estimate of $a$ in the following way.

If $n/4-1/n\leq a\leq n/4$, from \eqref{Aux-Term-Mono-1}, we have
\[
\begin{split}\frac{d}{da}\frac{(12a+n-2)^2}{128(n-2+4a)a}
&=\frac{12a+n-2}{128(n-2+4a)^2a^2}(n-2)(4a-n+2),\\
&\geq \frac{12a+n-2}{128(n-2+4a)^2a^2}(n-2)(2-\frac{4}{n})\\
&>0,
\end{split}
\]
 thus
 \[\frac{(12a+n-2)^2}{128(n-2+4a)a}\leq \frac{(12a+n-2)^2}{128(n-2+4a)a}\Big|_{a=\frac{n}{4}}=\frac{(2n-1)^2}{16n(n-1)}.\]
 Go back to  \eqref{Aux-a-n-cos-2-3}, we conclude that 
 \begin{equation}\label{Aux-a-n-cos-11-5}
  \eqref{Aux-a-n-cos-2-3}
 \geq\frac{1}{2}-\Big(\frac{1}{8}+\frac{(2n-1)^2}{16n(n-1)}\Big)\frac{(n-2)(n-4a)}{a}\geq 0,
 \end{equation}
  if \[a\geq \frac{ n(n-2)(6n^2-6n+1)}{4(3n-2)(2n^2-4n+1)}=
   \frac{n}{4}-\frac{n^2(n-1)}{2(3n-2)(2n^2-4n+1)}.\]
At this point, we would like point out that 
   \[\frac{n^2(n-1)}{2(3n-2)(2n^2-4n+1)}-\frac{1}{n}=\frac{(n-2)(n(n-1)(n-10)-2)}{2n(3n-2)(2n^2-4n+1)}\begin{cases}
      >0,&~\text{if}~n\geq 11, \\
     <0,&~\text{if}~4\leq n\leq 10.
\end{cases}\]

In summary,  by denoting 
\[\epsilon_n= \begin{cases}
     \frac{1}{n},&~\text{if}~n\geq 11, \\
    \frac{n^2(n-1)}{2(3n-2)(2n^2-4n+1)},&~\text{if}~4\leq n\leq 10,
\end{cases}\]
we conclude that  \eqref{Aux-a-n-sin-1}, \eqref{Aux-a-n-cos-2}  and  \eqref{Aux-a-n-sincos-3} holds simultaneously, and therefore \eqref{Diff-F-For-Q-neg-1} holds,  provide
 $n/4-\epsilon_n\leq a\leq n/4$. Since  $\overline{\Omega(a)}$  is   closed,  convex  and $O(n)$-invariant,  
 $\overline{\Omega(a)}$, for $n/4-\epsilon_n\leq a\leq n/4$, is preserved by Hamilton ODE  and further preserved  Ricci flow by Proposition \ref{RF-Inv-Tangent-Cone-1} and \ref{Hamilton'-maximum-principle}.
This is also true for  $\Omega(a)$  since positive scalar curvature is preserved  by Hamilton ODE  and   Ricci flow.  This finishes the proof of Theorem 
  \ref{Main-Theorem-1}.
 \end{proof}   
   
\subsection{Proof of Theorem \ref{Main-Convergence-Thm}}
 \begin{proof}[Proof of Theorem \ref{Main-Convergence-Thm}] 
 From  the definition of $\Omega(a)$ in   \eqref{B-W-Cone-Def-1},  
it is easy to see that  the family  $\overline{\Omega(a)}_{a\in [n/4-\epsilon_n,n/4)}\subset
S_B^2(\mathfrak{so}(n))$, where $\epsilon_n$ is defined in \eqref{BW-Conj-epsilon-n-Def},
satisfies 
\begin{enumerate}[(1)]
\item $\overline{\Omega(a)}$ is  closed convex $O(n)$-invariant cones  of full dimension,
\item  $a\mapsto \overline{\Omega(a)}$ is  continuous, 
 \item each $R\in \overline{\Omega(a)}\setminus \{0\}=\Omega(a)$ has positive scalar curvature,
\item $\overline{\Omega(a)}$ converges in the pointed Hausdorff topology to the one-dimensional cone $\mathbb{R}^+ I$ as $a\rightarrow n/4$.
\end{enumerate}
   Moreover,  when $n/4-\epsilon_n\leq a < n/4$,     the strict inequalities  in \eqref{Aux-a-n-cos-2-5} and \eqref{Aux-a-n-cos-11-5} hold, thus we have  that   the strict inequality in \eqref{Aux-a-n-cos-2} holds.  When $n/4-\epsilon_n< a < n/4$, then  the strict inequalities in \eqref{Aux-a-n-sin-eps-2} and  \eqref{Aux-a-n-sin-1} hold. Therefore,  in case of $n/4-\epsilon_n< a < n/4$, the strictly inequalities  in  \eqref{Diff-F-For-Q-nonnegative-5}  and \eqref{Diff-F-For-Q-neg-3} hold, i.e., 
  \begin{equation}\label{Diff-F-For-Q-neg-5}
dF_R(Q(R))< 0,\quad  \forall R\in \partial\Omega(a).
\end{equation}
In other words, we have 
 \begin{enumerate}[(1)]
 \setcounter{enumi}{4}
\item   \label{PF-TC-In-Cri-1} $Q(R)=R^2+R^\sharp$ is contained in the interior of the tangent cone of $\Omega(a)$ at $R$ for all $R \in \overline{\Omega(a)}\setminus \{0\}=\Omega(a)$ and all $a\in (n/4-\epsilon_n,n/4)$.
\end{enumerate}

Therefore, the  continuous family $\overline{\Omega(a)}_{a\in [n/4-\epsilon_n,n/4)}\subset
S_B^2(\mathfrak{so}(n))$ is a pinching family in the sense of Definition \ref{Pinching-family-Def}.  
By the general convergence Theorem \ref{Hamilton-BW-CP-RF},  we conclude that  on  a compact Riemannian manifold $(M,g)$ with the curvature operator at each point is  contained in     the interior of $\Omega(n/4-\epsilon_n)$,   the normalized Ricci flow  evolves $g$ to a constant curvature limit metric. Therefore, $M$ is  diffeomorphic to a spherical space form.  This finishes the proof of Theorem \ref{Main-Convergence-Thm}.
\end{proof}
   
  \section{K\"ahler manifolds and principal circle bundle}
As discussed in the introduction, one can investigate the topology  of K\"ahler manifolds through constructing  and studying the  principal circle bundle  over the  K\"ahler manifolds.  This    idea   tracks back  to   Kobayashi \cite{MR0154235}. 

\subsection{Geometry of principle  circle bundle}
Let us  assume  the existence of   principle circle  bundle $P$  over $(M^n,g)$ with a connection $\gamma$ for  the moment.
Following \cite{MR0154235}, for any $t>0$, we define a metric $\tilde{g}=\pi^*g+t^2\gamma^2$ on $P$, and
 express the curvature of $\tilde{g}$ in terms of those of $g$ and $\gamma$.  We shall also agree on that indices $i, j, k$ run from $1$ to $n$ and indices $\alpha, \beta, \lambda$ and $\mu$ run from $0$ to $n$.

Let $U$ be a open set in $M$ in which $g= \sum_{j=1}^n(\theta^j)^2.$ Let
 $(\omega_{ij})$ be a skew-symmetric matrix of $1$-forms which defines the Riemannian connection of $M$ on $U$ so that we have the following Cartan's structure equations:
\begin{equation}\label{Structure-Equ-1}
d\theta^i=-\sum_j\omega_{ij}\wedge \theta^j,
\end{equation}
\begin{equation}\label{Structure-Equ-2}
d\omega_{ij}=-\sum_k\omega_{ik}\wedge \omega_{kj}+\Omega_{ij},
\end{equation}
with 
\[\Omega_{ij}=\frac{1}{2}\sum_{k,l}K_{ijkl} \theta^k\wedge  \theta^l,\]
where $K_{ijkl}$ are the components of the curvature tensor with respect to $\theta^1,\cdots, \theta^n$. 

Since the structure group $S^1$ of $P$ is Abelian, the curvature of the connection $\gamma$ is defined as 
\begin{equation}\label{P-B-CC-Equ-2}
d\gamma=\pi^*A,
\end{equation}
where $A=\frac{1}{2} \sum_{i,j} A_{ij} \theta^i\wedge \theta^j,~A_{ij}=-A_{ij}$, is closed $2$-form on $M$. The cohomology class $[A]\in H^2(M,\mathbb{Z})$ is independent of the choice of connections and is called the characteristic class of $P$.

The first covariant derivative of $A_{ij}$ is defined by 
\begin{equation}\label{Cov-Der-Con-2}
\sum_kA_{ij,k}\theta^k=dA_{ij}-A_{kj}\omega_{ki}-A_{ik}\omega_{kj},
\end{equation}
and the second Bianchi identity asserts that 
\[A_{ij,k}+A_{jk,i}+A_{ki,j}=0.\]
For $\tilde{g}=\pi^*g+t^2\gamma^2$, set $\theta^0=t \gamma$,  then $\{\theta^\alpha\}_{\alpha=0}^n=\{t \gamma, \theta^1,\cdots, \theta^n\}$ is an orthonormal frame of $T^*P$. The Levi-Civita connection $1-$from $\tilde{\omega}_{\alpha\beta}$ is uniquely determined by
\begin{equation}\label{Structure-Equ-5}
d\theta^\alpha =-\tilde{\omega}_{ \alpha \beta }\wedge\theta^\beta, \quad \tilde{\omega}_{ \alpha \beta }=-\tilde{\omega}_{\beta  \alpha }.
\end{equation} 
Therefore, from \eqref{Structure-Equ-1}, \eqref{P-B-CC-Equ-2} and \eqref{Structure-Equ-5} we have, 
\begin{align*}
  -\sum_j\omega_{ij}\wedge \theta^j&= d\theta^i = -\tilde{\omega}_{ij}\wedge \theta^j-t\tilde{\omega}_{i0}\wedge \gamma,\\
t\big(\frac{1}{2}A_{ij} \theta^i\wedge\theta^j&\big)  =d\theta^0 =-\tilde{\omega}_{0 i}\wedge\theta^i,
\end{align*}
hence,
\begin{equation}\label{Connection-BP-Equ-1}
\tilde{\omega}_{ij}=\omega_{ij}-\frac{t^2}{2}A_{ij} \gamma, \quad \tilde{\omega}_{0 i}=\frac{1}{2}tA_{ij} \theta^j.
\end{equation} 
The Riemannian curvatures can be computed using the equation
\[d\tilde{\omega}_{ \alpha \beta }+\tilde{\omega}_{ \alpha \mu}\wedge\tilde{\omega}_{\mu \beta }=\frac{1}{2}R_{ \alpha \beta \mu \delta}\theta^\mu \wedge \theta^\delta.\]
Combine with \eqref{Structure-Equ-2}, \eqref{Cov-Der-Con-2} and \eqref{Connection-BP-Equ-1}, we have
\begin{align*}
   \frac{1}{2}R_{ij \alpha \beta }\theta^ \alpha \wedge \theta^\beta &=\frac{1}{2}K_{ijkl}\theta^k\wedge \theta^l-\frac{t^2}{4}(A_{ij}  A_{kl} +A_{ik}  A_{jl} )\theta^k\wedge \theta^l \\
    &\quad -\frac{t}{2}A_{ij,k} \theta^k\wedge \theta^0,\\
\frac{1}{2}R_{i0  \alpha \beta }\theta^\alpha \wedge \theta^\beta &=\frac{1}{4}(t^2A_{ik} A_{jk} )\theta^j\wedge \theta^0+\frac{1}{2}A_{ij,k}  \theta^j\wedge \theta^k.
\end{align*}
Therefore, we obtain 
\begin{proposition}\cite[Proposition 3]{MR0154235}\label{Curvature-PC-qJ-For-1}
The curvature of  $\tilde{g}$ on $P$ takes the form 
\begin{align}
   R_{ijkl} &=K_{ijkl}-\frac{t^2}{4}(2A_{ij}  A_{kl} +A_{ik}  A_{jl} -A_{il}  A_{jk} ),\\
   R_{ijk0} &=-\frac{t}{2}A_{ij,k},\\
   R_{i0 j0}&=\frac{t^2}{4}A_{ik} A_{jk}.
\end{align}
\end{proposition}

\subsection{Construction of principle circle bundle over a  K\"ahler manifold}  
Now, we turn to the construction of principle circle bundle over a  K\"ahler manifold $(M^m,g,J,\omega)$.  

Firstly, it was well-known that the set $P(M,S^1)$ of principle circle bundles over $M$  forms an additive group,  and the mapping which sends the principle circle bundles $P$
 to its characteristic class is an isomorphism of $P(M,S^1)$ onto  $H^2(M,\mathbb{Z})$ (e.g., see \cite{MR80919}). 
If $(M,g)$ has positive  bisectional curvature,  then   $H^2(M;\mathbb{R})\simeq \mathbb{R}$ (\cite[Theorem 4]{MR0227901}).   In this case,  
 the group $P(M,S^1)$ of all principal circle bundle  is isomorphic to the additive group of integers $\mathbb{Z}$.
In particular, for any  $q\in \mathbb{R}$ such that 
$[q\omega]\in H^{1,1}(M,\mathbb{C})\cap H^2(M,\mathbb{Z})$,  there  exists a  principal circle bundle $\pi: P\rightarrow M$ and a connection $\gamma$ on $P$ whose curvature form is $\pi^*(q\omega)$ (see \cite{MR80919}, and also \cite{MR112160}).  In the language of   K\"ahler  geometry,   $(M,q\omega)$ is a Hodge manifold and automatically projective algebraic, the K\"ahler class $q\omega$  defines an ample line bundle $L$ over $M$ such that the cone $P\times \mathbb{R}$ is identified with the complement of the zero section in $L^{-1}$, and  the total space  of the unit circle bundle in line bundle $L^{-1}$ is $P$ (see \cite{MR4761837}). For example, on $(\mathbb{P}^m,\omega_{\rm FS})$  (see \cite[Section 7]{MR80919}),  we have $P(M,S^1)\simeq H^2(\mathbb{P}^m;\mathbb{Z})=\mathbb{Z}$. 
The total space  of the unit circle bundle in line bundle $\mathcal{O}(-1)$ over $\mathbb{P}^m$ is $P= S^{2m+1}$,  which has the characteristic class represented by (curvature) $-\omega_{\rm FS}$;  while  the total space  of the unit circle bundle in the canonical line bundle   $\mathcal{O}(-(m+1))$  is $P= S^{2m+1}/\mathbb{Z}_{m+1}$.  

Let $\tau$ be any fixed positive number, consider the metric    $\tilde{g}=\pi^*g+(\frac{\tau}{q})^2\gamma^2$ on $P$.  In this case, 
denote $R$ to be  the Riemannian curvature operator   of  $\tilde{g}$  as in Proposition \ref{Curvature-PC-qJ-For-1} with  $A=q\omega$, i.e., $A_{ij} = J_{ij}$,  then $R$ have  the following relative simple from:
\begin{align}\label{Curvature-Omega-For-1}
   R_{ijkl} &=K_{ijkl}-\frac{\tau^2}{4}(2J_{ij}  J_{kl} +J_{ik} J_{jl} -J_{il}  J_{jk} ),\notag\\
   R_{ijk0} &=0,\\
   R_{i0 j0}&=\frac{\tau^2}{4}\delta_{ij}.\notag
\end{align}
Moreover, by  \eqref{Kahler-operator-CPn-1}, we may rewrite \eqref{Curvature-Omega-For-1} in the following form,
\begin{align}\label{Curvature-Omega-For-3}
  R =\frac{\tau^2}{4}I+(K-\frac{\tau^2}{2}E),
\end{align}
where we extend   the K\"ahler curvature $K$ on $P$  such that it is  evaluated to be zero whenever at least one of the vectors is tangent to the fiber, i.e.,  $K_{ijk0}=E_{ijk0} =0,
K_{i0 j0} =E_{i0 j0}=0.$ 

In this case, the Ricci curvature of $R$  in \eqref{Curvature-Omega-For-1}  is given by 
\begin{equation}\label{Ricci-Circle-Base-Relation-1}
\begin{split}
 R_{ij}&= K_{ij}-\frac{1}{2}\tau^2\delta_{ij}, \\
   R_{i0}&=0,\\
   R_{00}&=\frac{m}{2}\tau^2,
   \end{split}
\end{equation}
here $K_{ij}=\sum_{k=1}^{2m}K_{ikjk}$  is the Ricci curvature of $K$. 

The scalar curvature is given by 
\begin{equation}\label{Scalar-Circle-Base-Relation-1}
{\rm scal}(R)={\rm scal}(K)-\frac{m}{2}\tau^2,\quad \bar{\lambda}(R)=\frac{2m}{2m+1}(\bar{\lambda}-\frac{1}{4}\tau^2),
\end{equation}
here $\bar{\lambda}={\rm scal}(K)/(2m)$.

In general,   if there does not exist   $q\in \mathbb{R}$ such that 
$[q\omega]\in H^{1,1}(M,\mathbb{C})\cap H^2(M,\mathbb{Z})$,   one can  approximate  the  K\"ahler class   by  rational classes.
In fact, since $H^2(M;\mathbb{Z})$ form a basis in $H^2(M;\mathbb{R})$, 
 and the set  $\{bA;A\in H^2(M;\mathbb{Z}),b\in \mathbb{R}\}$ is dense in $H^2(M;\mathbb{R})$,  one can approximate $q\omega$ by some  integral class $A$ for some $q(=1/b)$.  Then there exists  principle circle bundle $P$ and a connection $\gamma$  such that the curvature  of $\gamma$  is  $\pi^*A$, and 
the curvature of $\tilde{g}=\pi^*g+(\tau/q)^2\gamma^2$ on $P$ is  arbitrarily close  to $R$  that is defined in \eqref{Curvature-Omega-For-1}. This is made precise in the following proposition.

\begin{proposition}\cite[Proposition 7-9]{MR0154235} \label{Stability-Curvature-Prop} On a  K\"ahler manifold $(M, g, J,\omega)$,  let $\tau$ be any fixed positive number, for any   positive number $\epsilon$, there exists a positive number $q$, and  a principle circle bundle $P$ over $M$ with a connection form $\gamma$,  such that
\begin{enumerate}
  \item the curvature of the connection $\gamma$  can be written as
   $d\gamma=\pi^*A$,
where $A=\frac{1}{2}A_{ij} \theta^i\wedge \theta^j$ is a harmonic form representing an element of  $H^2(M;\mathbb{Z})$ and satisfies
\[\text{maximum of}~\Big(\sum_{i,j}|A_{ij}-qJ_{ij}|^2+\sum_{i,j,k}|A_{ij,k}|^2\Big)<\frac{1}{q^2}\epsilon,\] 
  \item the curvature  $R'$ of  $\tilde{g}=\pi^*g+(\frac{\tau}{q})^2\gamma^2$ on $P$ satisfies
\[\text{maximum of}~\|R'- R\|<\epsilon.\]
\end{enumerate}
\end{proposition}

 \subsection{Algebraic propositions of the curvature norm}  
%
 Throughout this subsection,  let  $R$  be  the Riemannian curvature operator  which  is defined in \eqref{Curvature-Omega-For-1}. Moreover, set  $n=2m+1$.  Now we collect some formulas that express the norm of $R$ in terms of the norm of $K$. 
 \begin{proposition}
The norm of the Ricci  tensor of $R$ is
\begin{equation}\label{BP-Ricci-norm-1}
\begin{split}
|{\rm Ric}(R)|^2&=|{\rm Ric}(K)|^2-2m\tau^2\bar{\lambda}+ \frac{m^2+2m}{4}\tau^4.
\end{split}
\end{equation}
\end{proposition}
\begin{proof} From \eqref{Ricci-Circle-Base-Relation-1} we calculate
\begin{equation*} 
\begin{split}
|{\rm Ric}(R)|^2&=\sum_{\alpha,\beta}R_{\alpha\beta}^2=\sum_{i,j}R_{ij}^2+R_{00}^2+2\sum_{i}R_{i0}^2\\
&=\sum_{i,j}\big(K_{ij}-\frac{1}{2}\tau^2\delta_{ij}\big)^2+\big(\frac{m}{2}\tau^2\big)^2\\
&=|{\rm Ric}(K)|^2-2m\tau^2\bar{\lambda}+ \frac{m^2+2m}{4}\tau^4.
\end{split}
\end{equation*}
\end{proof}

\begin{proposition}
The norm of the trace-free Ricci tensor is 
\begin{equation}\label{BP-trace-free-Ricci-norm-1}
\begin{split}
|{\rm Ric}(R)_0|^2
&=|{\rm Ric}(K)_0|^2+\frac{2m}{2m+1}\Big(\bar{\lambda}-\frac{m+1}{2}\tau^2\Big)^2.
\end{split}
\end{equation}
\end{proposition}
\begin{proof} Plug  \eqref{Riemann-Ricci-Norm-1} and \eqref{Kaehler-Ric-Norm-1} into \eqref{BP-Ricci-norm-1},
we  obtain \eqref{BP-trace-free-Ricci-norm-1} from \eqref{Scalar-Circle-Base-Relation-1}.
\end{proof}

\begin{proposition}The norm of the  scalar part of   curvature  tensor $R$ is 
\begin{equation}\label{BP-Scalar-FullCurvature-norm-1}
\begin{split}
\|R_I\|^2
&=\frac{m}{2m+1}(\bar{\lambda}-\frac{1}{4}\tau^2)^2.
\end{split}
\end{equation}
\end{proposition}
\begin{proof} It is easy to see from  \eqref{curvature-norm-decomposition-3}  and \eqref{Scalar-Circle-Base-Relation-1}.
\end{proof}

\begin{proposition}
The norm of the  trace-free Ricci  part of  the curvature $R$ is  
\begin{equation}\label{BP-Ricci-FullCurvature-norm-1}
\begin{split}
\|R_{{\rm Ric}_0}\|^2=\frac{1}{2m-1}\Big(\frac{m+2}{2}\|K_{{\rm Ric}_0}\|^2+\frac{2m}{2m+1}\Big(\bar{\lambda}-\frac{m+1}{2}\tau^2\Big)\Big).\end{split}
\end{equation}
\end{proposition}
\begin{proof}
By \eqref{curvature-norm-decomposition-5},    \eqref{Ricci-Norm-RK-Rel-1} and \eqref{BP-trace-free-Ricci-norm-1}, we have 
\begin{equation*}
\|R_{{\rm Ric}_0}\|^2=\frac{1}{2m-1}|{\rm Ric}(R)_0|^2=\frac{1}{2m-1}\Big(2|{\rm Ric}(K)_0|_K^2+\frac{2m}{2m+1}\Big(\bar{\lambda}-\frac{m+1}{2}\tau^2\Big)\Big),
\end{equation*}
then  we obtain \eqref{BP-Ricci-FullCurvature-norm-1} from \eqref{Kaehler-Curvature-Ric-Norm-1}. 
\end{proof}

\begin{proposition}

The  norm of  the  curvature  tenor $R$ is  
\begin{equation}\label{BP-FullCurvature-norm-1}
\begin{split}
\|R\|^2&=\|K\|^2-\frac{3}{2}m\tau^2\bar{\lambda}+\frac{6m^2+5m}{16}\tau^4.
\end{split}
\end{equation}
\end{proposition}
\begin{proof} From \eqref{Norm-Full-Curvature-Def-1} we calculate
\begin{equation}\label{PB-FullCurvature-norm-1}
\|R\|^2 =\frac{1}{4}\sum_{\alpha,\beta,\gamma,\delta=0}^{2m}R_{\alpha\beta\gamma\delta}^2 =\frac{1}{4}\Big(4\sum_{i,j,k=1}^{2m}R_{ijk0}^2+4\sum_{i,j=1}^{2m}R_{0i0j}^2+\sum_{i,j,k,l=1}^{2m}R_{ijkl}^2\Big).
\end{equation}
From \eqref{Curvature-Omega-For-1}, we have
\begin{align}
\sum_{i,j,k=1}^{2m}R_{ijk0}^2&=0,\label{PB-Cur-Norm-Aux-1}\\
\sum_{i,j=1}^{2m}R_{0i0j}^2&=\sum_{i,j=1}^{2m}(\frac{\tau^2}{4}\delta_{ij})^2=\frac{m\tau^4}{8},\label{PB-Cur-Norm-Aux-3}\\
\sum_{i,j,k,l=1}^{2m}R_{ijkl}^2
&= \sum_{i,j,k,l=1}^{2m}\big(K_{ijkl}-\frac{\tau^2}{4}(2J_{ij}  J_{kl} +J_{ik} J_{jl} -J_{il}  J_{jk} )\big)^2\notag\\
&=\sum_{i,j,k,l=1}^{2m}K_{ijkl}^2-\tau^2(\sum_{i,k=1}^{2m}K(e_i, Je_i,e_k, Je_k)+\sum_{i,j=1}^{2m}K(e_i,e_j, e_i, e_j))\notag\\
&\quad +\frac{6m^2+3m}{4}\tau^4\notag\\
&=4\|K\|^2-3\tau^2{\rm scal}(K)+\frac{6m^2+3m}{4}\tau^4,\label{PB-Cur-Norm-Aux-5}
\end{align}
where we used the following identity in the last equality \eqref{PB-Cur-Norm-Aux-5},
\[\sum_{i,k=1}^{2m}K(e_i, Je_i,e_k, Je_k)=4\sum_{i,k=1}^{m}K(e_i, Je_i,e_k, Je_k)=2\,{\rm scal}(K).\]
Plugging \eqref{PB-Cur-Norm-Aux-1}, \eqref{PB-Cur-Norm-Aux-3} and \eqref{PB-Cur-Norm-Aux-5} into \eqref{PB-FullCurvature-norm-1}, we obtain \eqref{BP-FullCurvature-norm-1}.
\end{proof}

\begin{proposition}
The norm of the scalar flat part of   curvature  tensor $R$ is
\begin{equation}\label{BP-Scalar-flat-FullCurvature-norm-1}
\begin{split}
\|R_{{\rm Ric}_0}\|^2+\|W\|^2
&=\|K_{{\rm Ric}_0}\|^2+\|B\|^2 +\frac{m(3m+1)}{(m+1)(2m+1)}\big(\bar{\lambda}-\frac{m+1}{2}\tau^2\big)^2.\\
\end{split}
\end{equation}
\end{proposition}
\begin{proof} With \eqref{BP-Scalar-FullCurvature-norm-1} and \eqref{Kaehler-Curvature-Scalar-Norm-1}, we 
derive  from \eqref{BP-FullCurvature-norm-1} that 
\[\big(\|R\|^2-\|R_I\|^2\big)-\big(\|K\|^2-\|K_E\|^2\big)=\frac{m(3m+1)}{(m+1)(2m+1)}\big(\bar{\lambda}-\frac{m+1}{2}\tau^2\big)^2.\]
Thus we conclude \eqref{BP-Scalar-flat-FullCurvature-norm-1}.
\end{proof}
\begin{proposition}The norm of the Weyl curvature part of   curvature  tensor $R$ is
\begin{equation}\label{Weyl-Norm-PCB-1}
\|W\|^2= \|B\|^2+\frac{3m-4}{2(2m-1)}\|K_{{\rm Ric}_0}\|^2+ \frac{3m(m-1)}{(m+1)(2m-1)} \big(\bar{\lambda}-\frac{m+1}{2}\tau^2\big)^2.
\end{equation}
\end{proposition}
\begin{proof}From \eqref{BP-Scalar-flat-FullCurvature-norm-1} and \eqref{BP-Ricci-FullCurvature-norm-1}, we calculate 
\begin{equation*}
\begin{split}
\|W\|^2&=\|K_{{\rm Ric}_0}\|^2+\|B\|^2-\|R_{{\rm Ric}_0}\|^2\\
&=\|K_{{\rm Ric}_0}\|^2+\|B\|^2 +\frac{m(3m+1)}{(m+1)(2m+1)}\big(\bar{\lambda}-\frac{m+1}{2}\tau^2\big)^2\\
&\quad -\frac{1}{2m-1}\Big(\frac{m+2}{2}\|K_{{\rm Ric}_0}\|^2+\frac{2m}{2m+1}\big(\bar{\lambda}-\frac{m+1}{2}\tau^2\big)^2\Big)\\
&= \|B\|^2+\frac{3m-4}{2(2m-1)}\|K_{{\rm Ric}_0}\|^2+ \frac{3m(m-1)}{(m+1)(2m-1)} \big(\bar{\lambda}-\frac{m+1}{2}\tau^2\big)^2.
\end{split}
\end{equation*}
The proof is finished.
\end{proof}

\section{K\"ahler manifolds with  pinched curvature operators}
To make  the norm of the scalar flat part, $\|R_{{\rm Ric}_0}\|^2+\|W\|^2$,  of   curvature  tensor $R$ to be relatively small, according to  the formula \eqref{BP-Scalar-flat-FullCurvature-norm-1},  
 we may choose $\tau$ such that  
\begin{equation}\label{Minimal-Curvature-norm-a-1}
\frac{m+1}{2}\tau^2=\bar{\lambda}_0:=\frac{1}{{\rm Vol}(M)}\int_M\bar{\lambda}=\frac{1}{{2m\rm Vol}(M)}\int_M {\rm scal}(K).
\end{equation}
In this case, from \eqref{Scalar-Circle-Base-Relation-1},  we have
\begin{equation}\label{PC-Scalar-Relation-1}
\bar{\lambda}(R)=\frac{2m}{2m+1}(\bar{\lambda}-\frac{1}{2(m+1)}\bar{\lambda}_0).
\end{equation}
Furthermore, from  \eqref{BP-Scalar-FullCurvature-norm-1},  we have 
\begin{equation}\label{PC-Scalar-norm-1}
\begin{split}
\|R_I\|^2&=\frac{m}{2m+1}(\bar{\lambda}-\frac{1}{2(m+1)}\bar{\lambda}_0)^2\\
&=\frac{m(2m+1)}{4(m+1)^2}\Big(\bar{\lambda}^2+\frac{2}{2m+1}\bar{\lambda}(\bar{\lambda}-\bar{\lambda}_0)+\frac{1}{(2m+1)^2}(\bar{\lambda}-\bar{\lambda}_0)^2\Big).
\end{split}
\end{equation}
From \eqref{BP-Ricci-FullCurvature-norm-1}, we have
\begin{equation}\label{PC-Ricci-norm-1}
\begin{split}
\|R_{{\rm Ric}_0}\|^2
&=\frac{1}{2m-1}\Big(\frac{m+2}{2}\|K_{{\rm Ric}_0}\|^2+\frac{2m}{2m+1}\Big(\bar{\lambda}-\bar{\lambda}_0\Big)^2\Big).
\end{split}
\end{equation}
From \eqref{Weyl-Norm-PCB-1}, we have 
\begin{equation}\label{PC-Weyl-norm-1}
\|W\|^2= \|B\|^2+\frac{3m-4}{2(2m-1)}\|K_{{\rm Ric}_0}\|^2+ \frac{3m(m-1)}{(m+1)(2m-1)} \Big(\bar{\lambda}-\bar{\lambda}_0\Big)^2.
\end{equation}

\begin{proposition}\label{PC-Kahler-Pinching-Hui-1}  Assume $0<a<(2m+1)/4$.  Let  $K$ be an algebraic K\"ahler curvature operator, and let $R$ be the associated  algebraic Riemannian curvature operator    defined in \eqref{Curvature-Omega-For-1}. If we choose $\tau$ as specified in \eqref{Minimal-Curvature-norm-a-1}, then $R$ is contained in  the curvature cone $\Omega(a)$ (see \eqref{B-W-Cone-Def-1}) only if $K$ has positive scalar curvature and satisfies   the   pinching condition \eqref{Wilking-Kahler-Pinching-1}
for some $0<\epsilon<(2m+1)/2$. 
\end{proposition}
\begin{proof}
To ensure  that  $R$  is contained in  the curvature cone $\Omega(a)$, by \eqref{B-W-Cone-Def-1},  it is necessary for $R$ to have positive scalar curvature and satisfy the following  inequality, 
\begin{equation}\label{B-W-Kahler-Cone-1}
a \|R_{{\rm Ric}_0}\|^2+\Big(\frac{2m-1}{4}+a\Big)\|W\|^2\leq  \frac{2m+1-4a}{4}\|R_I\|^2.
\end{equation}
Using   the H\"older inequality, 
  \[|\bar{\lambda}(\bar{\lambda}-\bar{\lambda}_0)| \leq \frac{1}{4\epsilon}(\bar{\lambda}-\bar{\lambda}_0)^2+\epsilon\bar{\lambda}^2,~\forall \epsilon>0,\]
 we concude from  \eqref{PC-Scalar-norm-1} that 
 \begin{equation} 
\begin{split}
\|R_I\|^2&\geq \frac{m(2m+1)}{4(m+1)^2}\Big(\big(1-\frac{2}{2m+1}\epsilon\big)\bar{\lambda}^2+\frac{1}{(2m+1)^2}\big(1-\frac{2m+1}{2} \frac{1}{\epsilon}\big)(\bar{\lambda}-\bar{\lambda}_0)^2\Big).
\end{split}
\end{equation}
Therefore, combine with \eqref{PC-Ricci-norm-1} and \eqref{PC-Weyl-norm-1}, 
to ensure   the inequality \eqref{B-W-Kahler-Cone-1} holds,  it suffices to have the following inequality holds for some $0<\epsilon<(2m+1)/2$,
\begin{equation*}
\begin{split}
a & \frac{1}{2m-1}\Big(\frac{m+2}{2}\|K_{{\rm Ric}_0}\|^2+\frac{2m}{2m+1}\Big(\bar{\lambda}-\bar{\lambda}_0\Big)^2\Big)\\
&+\Big(\frac{2m-1}{4}+a\Big)\Big(\|B\|^2+\frac{3m-4}{2(2m-1)} \|K_{{\rm Ric}_0}\|^2+ \frac{3m(m-1)}{(m+1)(2m-1)} \Big(\bar{\lambda}-\bar{\lambda}_0\Big)^2\Big)\\
&\leq \frac{2m+1-4a}{4}
\frac{m(2m+1)}{4(m+1)^2}\Big(\big(1-\frac{2}{2m+1}\epsilon\big)\bar{\lambda}^2+\frac{1}{(2m+1)^2}\big(1-\frac{2m+1}{2} \frac{1}{\epsilon}\big)(\bar{\lambda}-\bar{\lambda}_0)^2\Big),
\end{split}
\end{equation*}
  i.e., $K$ satisfies the following inequality 
 \begin{equation}\label{B-W-Kahler-Cone-Aux-3}
 \begin{split}
 \Big(&\frac{3m-4}{8}+a\Big) \|K_{{\rm Ric}_0}\|^2+\Big(\frac{2m-1}{4}+a\Big)\|B\|^2\\
 &+\frac{m}{16(m+1)^2}\Big((12m^2-13)+4(6m+5)a +\frac{2 m+1 - 4 a }{2}\frac{1}{\epsilon}\Big) \Big(\bar{\lambda}-\bar{\lambda}_0\Big)^2\\
 &\leq \frac{2m+1-4a}{4} \frac{2m+1}{8(m+1)}\big(1-\frac{2}{2m+1}\epsilon\big) \|K_E\|^2,
  \end{split}
 \end{equation}
which is exactly  the pinching condition \eqref{Wilking-Kahler-Pinching-1}.

Furthermore, if $K$ satisfies \eqref{B-W-Kahler-Cone-Aux-3}, we have 
\begin{equation}\label{B-W-Kahler-Cone-Aux-5}
\begin{split}
\Big((12m^2-13)+4(6m+5)a\Big) \Big(\bar{\lambda}-\bar{\lambda}_0\Big)^2\leq (2m+1-4a)(2m+1)\bar{\lambda}^2,
\end{split}
\end{equation}
which implies that
\begin{equation}\label{Pos-Sca-Aux-1}
\frac{1}{(2m+1)^2}(\bar{\lambda}-\bar{\lambda}_0)^2 \leq \frac{2m+1-4a}{(2m+1)((12m^2-13)+4(6m+5)a)}\bar{\lambda}^2<\bar{\lambda}^2.
\end{equation}
 Additionally,  the equation \eqref{PC-Scalar-Relation-1}  can be rewritten equivalently as  
\begin{equation}\label{Pos-Sca-Aux-3} 
\bar{\lambda}(R)
=\frac{m}{m+1}\big(\bar{\lambda}+\frac{1}{2m+1}(\bar{\lambda}-\bar{\lambda}_0)\big).
\end{equation}
Consequently, if  $K$ has positive scalar curvature, i.e., $\bar{\lambda}>0$, we conclude that  $R$ has positive scalar curvature from \eqref{Pos-Sca-Aux-1} and \eqref{Pos-Sca-Aux-3},  i.e., $\bar{\lambda}(R)>0$. Therefore, if  $K$ has positive scalar curvature and satisfies  the pinching condition \eqref{Wilking-Kahler-Pinching-1} for some $0<\epsilon<(2m+1)/2$, then  $R$ is contained in  the curvature cone $\Omega(a)$. 
\end{proof}  

\begin{proof}[{Proof of  Theorem \ref{Kahler-Pinching-Bundle-3}}]

Now  let $(M,g,J,\omega)$ be a compact  K\"ahler manifold of complex dimension $m$.
Assume that the scalar curvature  is positive, i.e., $\bar{\lambda}>0$,  and the  K\"ahler curvature $K(p)$  at each point $p\in M$ satisfies the pinching condition
\eqref{Wilking-Kahler-Pinching-1} for some $0<\epsilon<(2m+1)/2$ and $(2m+1)/4-\epsilon_{2m+1}< a \leq (2m+1)/4$,   where $\epsilon_{2m+1}$ is defined in \eqref{BW-Conj-epsilon-n-Def}.  Let $R$ be the associated  algebraic Riemannian curvature operator which  is defined in \eqref{Curvature-Omega-For-1}.  Since $R$ depends only on $p\in M$ not on $\pi^{-1}(p)\subset P$, will denote $R$ by $R(p)$. If we choose $\tau$ as in \eqref{Minimal-Curvature-norm-a-1}, by Proposition  \ref{PC-Kahler-Pinching-Hui-1}, we conclude that  $R(p)$  at each point $p\in M$ is contained in  the curvature cone $\Omega(a)$ for $(2m+1)/4-\epsilon_{2m+1}< a \leq (2m+1)/4$.  

Furthermore, since $M$ is compact,  and $(2m+1)/4-\epsilon_{2m+1}< a \leq (2m+1)/4$, there is a some positive
number $\epsilon$ such that the $\epsilon-$neighbourhood of the  curvature set \[\{R(p), \text{the algebraic curvature operator   defined in \eqref{Curvature-Omega-For-1}},   p\in M\},\] is contained in  the  interior of the curvature cone $\Omega((2m+1)/4-\epsilon_{2m+1})$.

For such $\epsilon>0$ and $\tau$ choosen in \eqref{Minimal-Curvature-norm-a-1},  by  Proposition \ref{Stability-Curvature-Prop}, there exists a positive number $q$, and  a principle circle bundle $P$ over $M$ with a connection form $\gamma$,  such that the curvature  $R'$ of  $\tilde{g}=\pi^*g+(\frac{\tau}{q})^2\gamma^2$ on $P$ satisfy
\[\max \|R'-  R\|<\epsilon.\]
In particular, $R'$ is contained   in the interior curvature cone  $\Omega((2m+1)/4-\epsilon_{2m+1})$. By  the Differential Sphere Theorem \ref{Main-Convergence-Thm},  we conclude that the principle circle bundle $P$  is  diffeomorphic to a spherical space form. In particular,  the universal covering space of   $P$  is  diffeomorphic with $S^{2m+1}$.  As in Kobayashi \cite{MR0154235}, 
the   exact homotopy sequences of the fiber bundle  $S^1\rightarrow S^{2m+1} \rightarrow  \mathbb{P}^m$  and $S^1\rightarrow P \rightarrow  M$ give
an isomorphism  $\pi_i(M)=\pi_i(\mathbb{P}^m)$ for $i\geq 2$.   
Moreover, in case of $(2m+1)/4-\epsilon_{2m+1}< a \leq (2m+1)/4$, from Corollary \ref{Ricci-Positive-Cond-1}, we know that 
$\Omega(a)$ is contained  in the curvature cone with positive Ricci curvature.   
Furthermore,  for   $X,Y$ horizontal in the submersion $\pi: P \rightarrow  M$, we see from \eqref{Ricci-Circle-Base-Relation-1}
that \begin{equation*}\label{Ricci-Base-Bundle-Positive-1}
{\rm Ric}(K)(X,Y)={\rm Ric}(R)(X,Y)+\frac{1}{2}\tau^2g(X, Y)\geq \frac{1}{2}\tau^2 g(X, Y).
\end{equation*}
In particular, $M$ has positive Ricci curvature and thus is simply connected  by a theorem of Kobayashi \cite{MR0133086}. Since   $\mathbb{P}^m$ is also simply connected, we conclude that $\pi_i(M)\simeq \pi_i(\mathbb{P}^m)$ for all $i$. 
The  proof of Theorem  \ref{Kahler-Pinching-Bundle-3} is finished.
\end{proof}

\subsection*{Acknowledgments}
Part of this work was supported by NSFC grants No.11501285 and No.11871265. The author would like to thank Professor Yalong Shi for the useful discussions. 
\bibliographystyle{abbrvnat}

\end{document}